\documentclass[12pt]{article}
\usepackage{amsmath}
\usepackage{amssymb}
\usepackage{amscd}
\usepackage{amsthm}
\usepackage{mysects}

\def\no{\noindent}
\def\bul{\bullet}
\numberwithin{equation}{subsection}

\theoremstyle{plain}
\newtheorem{theorem}[equation]{Theorem}

\newtheorem{proposition}[equation]{Proposition}
\newtheorem{lemma}[equation]{Lemma}
\newtheorem{corollary}[equation]{Corollary}

\theoremstyle{remark}
\newtheorem{remark}[equation]{Remark}

\theoremstyle{definition}
\newtheorem{definition}[equation]{Definition}

\newcommand{\lra}{\longrightarrow}
\newcommand{\ra}{\rightarrow}
\newcommand{\restr}{\mbox{\Large \(|\)\normalsize}}
\newcommand{\N}{\mathbb N}
\newcommand{\R}{\mathbb R}

\newcommand{\Z}{\mathbb Z}

\newcommand{\acts}{\curvearrowright}

\def\D{\partial}
\newcommand{\al}{\alpha}
\def\de{\delta}
\def\De{\Delta}
\def\nbd{neighborhood}

\def\eps{\epsilon}
\def\ga{\gamma}
\def\Ga{\Gamma}

\def\lra{\longrightarrow}

\def\larrow{\leftarrow}

\newcommand{\oa}{\overrightarrow}

\def\om{\omega}

\def\ra{\rightarrow}

\def\si{\sigma}
\def\Si{\Sigma}

\def\be{\beta}

\def\th{\theta}

\def\geo{\partial_{\infty}}
\def\t{\tilde}
\def\defeq{:=}
\def\hook{\hookrightarrow}

\newcommand{\ol}{\overline}
\newcommand{\<}{\langle}
\renewcommand{\>}{\rangle}
\def\RA{\joinrel\relbar\joinrel\relbar\joinrel\relbar\joinrel\relbar\joinrel\relbar\joinrel\rightarrow}

\newcommand{\pres}[1]{0\larrow \Z\larrow P_0\larrow\ldots\larrow P_{#1}\larrow 0}

\begin{document}

\title{Coarse Alexander duality and duality groups}
\author{Misha Kapovich\thanks{Supported by NSF grant DMS-96-26633.}
\\
Bruce Kleiner\thanks{Supported by a Sloan Foundation
Fellowship, and NSF grants
DMS-95-05175, DMS-96-26911.}}
\date{November 1, 1999}
\maketitle

\begin{abstract}
We study discrete group actions on coarse Poincare duality spaces, e.g.
acyclic simplicial complexes which admit free cocompact group
actions by Poincare duality groups.  When $G$ is an $(n-1)$ dimensional
duality group and $X$ is a coarse Poincare duality space of formal
dimension $n$, then a free simplicial action $G\acts X$ determines
a collection of ``peripheral'' subgroups $F_1,\ldots, F_k\subset G$ so that the
group pair $(G;F_1,\ldots,F_k)$ is an $n$-dimensional Poincare duality pair.
In particular, if $G$ is a $2$-dimensional 1-ended group of type $FP_2$, 
and $G\acts X$ is a free simplicial action on a coarse $PD(3)$ space $X$,
then $G$ contains surface subgroups; if in addition $X$ 
 is  simply connected, then we obtain a partial generalization of the
Scott/Shalen  compact core theorem to the setting of coarse $PD(3)$
spaces.
\end{abstract}

\setcounter{section}{1}
\setcounter{subsection}{0}

\subsection{Introduction}

In this paper we study simplicial complexes which
behave homologically (in the large-scale) like
$\R^n$, and discrete group actions on them.
Our main objective is a partial generalization of the
Scott/Shalen compact core theorem for $3$-manifolds
(\cite{scott}, see also \cite{jaco})
to the setting of Poincare duality spaces of arbitrary
dimension.  In the one ended case, the 
compact core theorem 
says that if $X$ is a contractible $3$-manifold and $G$ 
 is a
finitely generated one-ended group acting discretely and freely
on $X$, then the
quotient $X/G$ contains a compact core -- a compact
submanifold with (aspherical) incompressible 
boundary $Q\subset X/G$ so that
the inclusion $Q\ra X/G$ is a homotopy equivalence.
The proof of the compact core theorem relies on standard
tools in $3$-manifold theory like transversality, which has 
no appropriate analog in the $3$-dimensional coarse Poincare duality
space setting, and the Loop Theorem, which has no analog 
even for manifolds when the dimension is at least $4$.

We now formulate our analog of the core theorem. 
For our purpose, the appropriate substitute for a finitely generated,
one-ended, $2$-dimensional group will be a {\em duality group}
of dimension $n-1$.   We recall \cite{Bieri-Eckmann1}
that a group $G$ is a $k$-dimensional 
duality group if $G$ is of type $FP$,
$H^i(G;\Z G)=0$ for $i\neq k$, and $H^k(G;\Z G)$ is torsion-free~\footnote{We never make use of the last assumption 
about $H^k(G;\Z G)$  in our paper.}.   Examples of duality groups include:

A. Freely indecomposable $2$-dimensional groups of type $FP_2$;
for instance, torsion free one-ended 1-relator groups.

B. The fundamental groups of compact aspherical manifolds with 
incompressible aspherical boundary \cite{Bieri-Eckmann1}. 

C. The product of two duality groups. 

D. Torsion free $S$-arithmetic groups \cite{BS1}.

\no
Instead of $3$-dimensional contractible manifolds, we work with
 a class of metric simplicial
complexes which we call ``coarse $PD(n)$
spaces''.   These will be defined in section
\ref{cpd}, but we note that important examples
include  universal covers of closed
aspherical $n$-dimensional PL-manifolds, 
acyclic complexes $X$ with $H^*_c(X)\simeq H^*_c(\R^n)$
which admit free cocompact simplicial group actions, and uniformly 
acyclic $n$-dimensional  PL-manifolds with bounded geometry.  
We recall that an $n$-dimensional Poincare duality group ($PD(n)$ group) is
a duality group $G$ with $H^n(G;\Z G)\simeq \Z$.  Our  group-theoretic
analog for the compact core will an $n$-dimensional
Poincare duality pair ($PD(n)$ pair), i.e. a group pair $(G;F_1,\ldots, F_k)$ whose
double with respect to the $F_i$'s is an $n$-dimensional Poincare
duality group, \cite{dunwoodydicks}.  In this case the ``peripheral''
subgroups $F_i$ are $PD(n-1)$ groups.
See section \ref{groupprelim} for more details.

\begin{theorem}
\label{mainduality}
Let $X$ be a coarse $PD(n)$ space, and let $G$ be an
$(n-1)$-dimensional duality group acting 
discretely and simplicially on $X$.  Then:

1.  $G$ contains subgroups $F_1,\ldots F_k$
(which are canonically defined up to conjugacy
by the action $G\acts X$) so that $(G;\{F_i\})$
is a  $PD(n)$  pair. 

2. There is a connected $G$-invariant subcomplex $K\subset X$
so that $K/G$ is compact, the stabilizer of each component
of $X-K$ is conjugate to one of the $F_i$'s, and each
component of $\ol{X-K}/G$ is one-ended. 
\end{theorem}

Thus, the duality groups $G$ which appear in the above theorem 
behave homologically like the groups in example B. 
As far as we know, Theorem \ref{mainduality} is new
even in the case that $X\simeq \R^n$, when $n\geq 4$.
Theorem \ref{mainduality} and Lemma \ref{thesearepdn}  imply

\begin{corollary}
Suppose a $PD(n)$ group $\Ga$ acts freely, simplicially, and cocompactly
on an acyclic complex.  Then any $(n-1)$-dimensional duality subgroup
$G\subset\Ga$ contains a finite collection $H_1,\ldots,H_k$ of $PD(n-1)$
subgroups so that the group pair $(G,\{H_i\})$ is a $PD(n)$ pair.
\end{corollary}
In the next theorem we obtain a partial generalization 
  of the  Scott-Shalen theorem for groups acting
on coarse $PD(3)$ spaces.

\begin{theorem}
\label{main3d}
Let $G\acts X$ be a free simplicial action of a $2$-dimensional,
one-ended group of type $FP_2$ on a simply connected 
coarse $PD(3)$ space $X$.
 Then there exists
a complex $Y$ and a proper homotopy equivalence
$f:X/G\ra Y$ which is a homeomorphism away from
a compact subset, where $Y=Q\cup (E_1\sqcup\ldots\sqcup E_k)$, and

1. $Q$ is a finite subcomplex of $Y$, and $Q\hook Y$ is a homotopy equivalence.

2. The $E_i$'s are disjoint and one-ended. 
For each $i$, $S_i\defeq E_i\cap Q$ is a closed aspherical surface,
and $S_i\hook E_i$ is a homotopy equivalence.

3. Each inclusion $S_i\hook Q$ is $\pi_1$-injective.

4. $(Q;S_1,\ldots,S_k)$ is a 
Poincare pair \cite{weinberger}. 
In particular, $Q$ is a finite  Eilenberg-MacLane space for $G$.
\end{theorem}

\begin{corollary}
If $G$ is a group of type $FP_2$, $dim(G)\leq 2$,  and $G$ acts
freely simplicially on a coarse $PD(3)$ space, then
either $G$ contains a surface group, or $G$ is free.
\end{corollary}
\proof
Let $G=F*(*_i G_i)$ be a free product decomposition
where $F$ is a finitely generated free group, and each
$G_i$ is finitely generated, freely indecomposable, and 
non-cyclic.  Then by Stallings' theorem on ends of groups,
each $G_i$ is one-ended, and hence is a $2$-dimensional
duality group.  By Theorem \ref{main3d}, each $G_i$ contains
surface subgroups.
\qed

\medskip
We believe that Theorem \ref{main3d} still holds if one 
relaxes the $FP_2$ assumption to finite generation, and
 we conjecture that any finitely generated 
group which acts freely, simplicially, but not
cocompactly, on a coarse $PD(3)$ space is finitely presented. 
We note that Bestvina and Brady \cite{bestbrad} construct
$2$-dimensional groups which are $FP_2$ but not finitely presented.

In Proposition \ref{peripheraluniqueness} we prove an analog 
of the uniqueness theorem for peripheral structure \cite{jacoshalen,johannson} 
for fundamental groups of acylindrical $3$-manifolds with
aspherical incompressible boundary.

We were led to Theorems \ref{mainduality}, \ref{main3d} by our 
earlier work on hyperbolic groups with one-dimensional boundary \cite{KK}; 
in that paper we conjectured that every torsion-free hyperbolic group $G$  
whose boundary is homeomorphic to the Sierpinski carpet is the fundamental group of 
a compact hyperbolic 3-manifold with totally geodesic boundary. 
In the same paper we showed that such a group $G$ is part of  a canonically 
defined $PD(3)$ pair and that our conjecture would follow if one 
knew that $G$ were a 3-manifold group. One approach to proving this  
is to produce an algebraic counterpart to the Haken hierarchy for Haken 
3-manifolds in the context of $PD(3)$ pairs. 
We say that a $PD(3)$
pair $(G;H_1,\ldots,H_k)$ is {\em Haken} if it admits
a nontrivial splitting\footnote{If $k>0$ then
such a splitting always exists.}. One would like to show that Haken 
 $PD(3)$ pairs always admit nontrivial splittings over $PD(2)$
pairs  whose peripheral structure is compatible with that of $G$. 
Given this, one can create a hierarchical decomposition of the 
group $G$, and try to show that the terminal groups correspond to
fundamental groups of $3$-manifolds with boundary.  The corresponding
$3$-manifolds might then be glued together along boundary surfaces
to yield a $3$-manifold with fundamental group $G$.
At the moment, the biggest obstacle in this hierarchy program 
appears to be the
first step; and the two theorems above
provide a step toward overcoming it under the
assumption that finitely generated subgroups of
$PD(3)$ groups are of type $FP_2$.

As an application of Theorems \ref{mainduality} and \ref{main3d} 
and the techniques used in their proof,
we give examples of $(n-1)$-dimensional groups which cannot
act discretely and simplicially on coarse $PD(n)$
 spaces (see section \ref{appl} for details):

1. A $2$-dimensional one-ended group of type $FP_2$ with
positive Euler characteristic cannot act on a coarse $PD(3)$
space.   For example,  the semi-direct
product of two finitely generated free groups.

2. For $i=1,...,\ell$ let $G_i$ be a duality group of dimension $n_i$ 
and assume that for $i=1,2$ the group $G_i$ is  not a $PD(n_i)$ group. Then the product 
$G_1 \times ... \times G_\ell$ cannot act on a coarse $PD(n)$
space where $n-1=n_1 + ... + n_\ell$.  The case when $n=3$ 
is due to Kropholler, \cite{kropholler}.


3.  If $G_1$ is a $k$-dimensional duality group and $G_2$ is the the Baumslag-Solitar
group $BS(p,q)$ (where $p\ne \pm q$), then 
the direct product $G_1\times G_2$ cannot  act on a coarse $PD(3+k)$ space.
 In particular, $BS(p,q)$ cannot act on 
a coarse $PD(3)$ space.

4.  An $(n-1)$-dimensional  group $G$ of type $FP_{n-1}$  
which contains infinitely many conjugacy classes of coarsely 
non-separating maximal $PD(n-1)$ subgroups. 

Our theme is  related to the problem of finding an
$n$-thickening of an  aspherical polyhedron $P$ up to homotopy,
i.e. finding a homotopy equivalence $P\to M$ where $M$
is a compact manifold with boundary and $dim(M)=n$.   If $k=dim(P)$
then we may immerse $P$ in $\R^{2k}$ by  general position,
and obtain a $2k$-manifold thickening $M$ by 
``pulling back'' a regular neighborhood.  Given  an $n$-thickening
$P\ra M$  we may construct a free simplicial action of $G=\pi_1(P)$
on a coarse   $PD(n)$ space by modifying the geometry of $Int(M)$
and passing to the universal cover. In particular, if $G$ cannot
act on a coarse $PD(n)$ space then no such $n$-thickening can 
exist.   In a subsequent paper with M.~Bestvina \cite{withmladen}  
we give examples of finite $k$-dimensional
aspherical polyhedra $P$ whose fundamental groups cannot act freely
simplicially on any coarse $PD(n)$ space for $n<2k$, and hence the
polyhedra $P$ do not admit $n$-thickening for $n<2k$.

In this paper we develop and use ideas in coarse topology which originated
in earlier work by a number of authors: 
\cite{blockwein,farbschwartz,gromov,higsonroe,roe,schwartz1,schwartz2}.
Other recent papers involving similar ideas include 
\cite{vaisala,farbeskin,farbmosher}, and especially \cite{bowditch},
which has considerable overlap with this paper.  An adaptation of Richard Schwartz's
coarse Alexander duality to coarse $PD(n)$ spaces plays an important role
in the proofs of our main results.

To give an idea of the proof of Theorem \ref{mainduality}, consider the case 
when the coarse $PD(n)$-space $X$ happens to be $\R^n$ with a 
uniformly acyclic bounded geometry
triangulation. We take combinatorial
tubular neighborhoods $N_R(K)$ of 
a $G$-orbit $K$ in $X$ and analyze the structure of connected components of
$X-N_R(K)$. Following R.~Schwartz  we call a connected component $C$ of 
$X-N_R(K)$ {\em deep} if $C$ is not contained in any  tubular neighborhood of 
$K$.  When $G$ is a group of type $FP_n$, 
using Alexander duality one shows that deep components of $X-N_R(K)$ 
{\em stabilize}: there exists $R_0$ so that no 
deep component of $X- N_{R_0}(K)$ breaks up into multiple
deep components as $R$ increases beyond $R_0$. 
If $G$ is an $(n-1)$-dimensional duality group then the idea is 
to show that the stabilizers of 
of deep components of $X- N_{R_0}(K)$ are $PD(n-1)$-groups, which 
is the heart of the proof. These groups define the peripheral subgroups 
$F_1,\ldots,F_k$ of the $PD(n)$ pair structure $(G;F_1,\ldots,F_k)$ for $G$.

When $X$ is a coarse $PD(n)$-space rather than $\R^n$,
 one does not have Alexander duality since 
Poincare duality need not hold locally. However there is a coarse version of 
Poincare duality which we use  to derive an appropriate coarse analogue of  
Alexander duality. Roughly speaking this goes as follows.  If $K\subset\R^n$
is a subcomplex then Poincare duality gives an
isomorphism
$$
H^*_c(K)\ra H_{n-*}(\R^n,\R^n-K).
$$
This fails when we replace $\R^n$ by a general coarse $PD(n)$
space $X$. We prove however that for a certain
constant $D$ there are homomorphisms defined on tubular neighborhoods of $K$:
$$
P_{R+D}: H^k_c(N_{D+R}(K))\to H_{n-k}(X, Y_R), \hbox{~where~} Y_R:=\ol{X-N_R(K)}, 
$$
which determine an {\em approximate isomorphism}.
This means that for every $R$ there is an $R'$ (one may take 
$R'=R+2D$) so that the homorphisms $a$ and $b$ in the following
commutative diagram are zero:
$$
\begin{array}{cccccc}
ker(P_{R'}) & \to & H_c^k(N_{R'}(K)) & \stackrel{P_{R'}}{\lra}  H_{n-k}(X, Y_{R' -D}) & \to & coker(P_{R'})\\
a\downarrow & ~ & \downarrow & \downarrow & ~ & b\downarrow\\
ker(P_{R}) & \to & H_c^k(N_{R}(K)) & \stackrel{P_{R}}{\lra}  H_{n-k}(X, Y_{R -D}) & \to & coker(P_{R})
\end{array}
$$
This  coarse version of
Poincare duality leads to coarse Alexander duality, which suffices for our purposes. 

\medskip
\no
{\bf Organization of the paper.}  In section 2 we introduce metric simplicial
complexes and recall notions from coarse
topology.    Section 3 reviews some facts and definitions from cohomological group
theory, duality groups, and group pairs.  In section 4 we define approximate
isomorphisms between inverse and direct systems of abelian groups,
and compare these with Grothendieck's pro-morphisms.  Section 5
provides finiteness criteria for groups, and 
establishes approximate isomorphisms between group cohomology
and cohomologies of nested families of simplicial complexes.   In section 
6 we define coarse $PD(n)$ spaces, give examples, and prove coarse
Poincare duality for coarse $PD(n) $ spaces.   In section 7 we prove
coarse Alexander   duality and apply it to coarse separation.
In section 8 we prove Theorems \ref{mainduality}, \ref{main3d},  Proposition
\ref{peripheraluniqueness}, and variants of Theorem \ref{mainduality}.
In section 9 we apply coarse Alexander duality and Theorem \ref{mainduality}
to show that certain groups cannot act freely simplicially on coarse 
$PD(n)$ spaces. In the Appendix (section 10) we give a brief account of coarse 
Alexander duality for uniformly acyclic triangulated manifolds of 
bounded geometry. The reader interested in manifolds and not in Poincare comlexes 
can use this as a replacement of Theorem \ref{aduality}.

\medskip
\no
{\bf Suggestions to the reader.} 
Readers familiar with Grothendieck's 
pro-morphisms may wish to read the second part of section 4, which
will allow them to translate statements about approximate isomorphisms
into pro-language.  Readers who are not already familiar
with pro-morphisms may simply skip this.   Those who are interested in 
finiteness properties of groups may find section 5, especially Theorems
\ref{finitenesstheorem} and Corollary \ref{finitenesscor}, of independent interest.

\medskip
\no
{\bf Acknowledgements.} We are grateful for M.~Bestvina and S.~Weinberger  
for useful conversations about coarse Poincare duality.

\tableofcontents

\subsection{Geometric Preliminaries}
\label{geomprelim}

{\bf Metric simplicial complexes.}
Let $X$ be the geometric realization of a connected
locally finite  simplicial complex. Henceforth we 
will conflate simplicial complexes with their geometric realizations.
 We will metrize the 1-skeleton $X^1$ 
of  $X$ by declaring each edge to have unit length and taking the 
corresponding path-metric. Such an
$X$ with the metric on $X^1$ will be called a {\em metric simplicial complex}. 
The complex $X$ is said to have {\em bounded geometry} if all links 
have a uniformly bounded number of simplices; this is equivalent to saying that the metric space $X^1$ 
is locally compact and every $R$-ball in $X^1$ can be covered by at most $C=C(R,r)$ $r$-balls 
for any $r>0$. In particular, $dim(X)<\infty$. If $K\subset X$ is a subcomplex  and $r$ is a positive integer 
then we define (combinatorial) $r$-tubular \nbd\ 
$N_r(K)$ of $K$ to be $r$-fold iterated closed star of $K$, $St^r(K)$; 
we declare $N_0(K)$ to be $K$ itself.   Note that for $r>0$, $N_r(K)$
is the closure of its interior.  The diameter of $K$ is defined to be 
the diameter of its zero-skeleton, and $\D K$ denotes the frontier of $K$, 
which is a subcomplex. For each vertex $x\in X$ and $R\in \Z_+$ we let 
$B(x,R)$ denote $N_R(\{x\})$, the ``$R$-ball centered at $x$''. 

\medskip
\no
{\bf Coarse Lipschitz and uniformly proper maps.} We recall that a map 
$f:X\to Y$ between metric spaces is 
called {\em $(L,A)$-Lipschitz} if 
$$
d(f(x), f(x'))\le L d(x,x') +A
$$
for any $x,x'\in X$. A map is {\em coarse Lipschitz} if it is $(L,A)$-Lipschitz for some $L, A$. A coarse Lipschitz map 
$f: X\to Y$ is called {\em uniformly proper} if there is a proper function $\phi: \R_+\to\R_+$ 
(a {\em distortion function}) such that
$$
d(f(x), f(x'))\ge \phi(d(x,x'))
$$
for all $x,x'\in X$. 

Throughout the paper we will use
simplicial (co)chain complexes and integer coefficients.
If $C_*(X)$ is the simplicial chain complex and $A\subset C_*(X)$, 
then the {\em support of  $A$}, denoted $Support(A)$,   is the smallest 
subcomplex $K\subset X$ so that $A\subset C_*(K)$. Throughout the paper 
we will assume that morphisms between simplicial chain complexes preserve 
the usual augmentation. 

If $X,Y$ are metric simplicial complexes as above then a homomorphism
$$
h: C_*(X)\to C_*(Y)
$$ 
is said to be {\em coarse Lipschitz} if for each simplex $\sigma\subset X$,
$Support(h(C_*(\sigma)))$ has uniformly bounded diameter. 
The {\em Lipschitz constant of $h$} is  
$$
\max_{\si} diam(Support(h(C_*(\sigma)))).
$$ 
A homomorphism $h$ is said to be {\em uniformly proper} if 
it is coarse 
Lipschitz and there exists a proper function $\phi: \R_+ \to \R_+$
(a {\em distortion function}) such that 
for each subcomplex $K\subset X$ of diameter $\ge r$,  $Support(h(C_*(K)))$ 
has diameter $\ge \phi(r)$. 
We will apply this definition only to chain mappings and chain homotopies\footnote{Recall that
there is a standard way to triangulate the product $\De^k\times [0,1]$; we can use this
to triangulate $X\times[0,1]$ and hence view it as a metric simplicial complex.}.
We say that a homomorphism $h:C_*(X)\ra C_*(X)$ has displacement $\leq D$
if for every simplex $\si\subset X$,  $Support(h(C_*(\si)))\subset N_D(\si)$. 

We may adapt all of the definitions from the previous paragraph to mappings
between other (co)chain complexes associated with metric simplicial complexes,
such as the compactly supported cochain complex $C^*_c(X)$.

\medskip
\no
{\bf Coarse topology.} A metric simplicial complex $X$ is said to be {\em uniformly acyclic} if 
for every $R_1$ there is an $R_2$ such that for each subcomplex $K\subset X$ of diameter $\le R_1$ the 
inclusion $K\to N_{R_2}(K)$ induces zero on  reduced homology groups. 
Such a  function $R_2=R_2(R_1)$
will be called an {\em acyclicity function} for $C_*(X)$.
Let $C^*_c(X)$ denote the complex of simplicial cochains, and suppose
$\al:C^n_c(X)\ra\Z$ is an augmentation for $C^*_c(X)$.  Then the pair
$(C^*_c(X),\al)$ is {\em uniformly acyclic} if there is
an $R_0>0$ and a function $R_2=R_2(R_1)$ so that
for all $x\in X^0$ and all $R_1\geq R_0$,
$$Im(H^*_c(X,\ol{X-B(x,R_1}))\ra H^*_c(X,\ol{X-B(x,R_2)}))$$
maps isomorphically onto $H^*_c(X)$ under 
$H^*_c(X,\ol{X-B(x,R_2)})\ra H^*_c(X)$, and $\al$ 
induces an isomorphism $\bar\al:H^n_c(X)\ra\Z$.

Let $K\subset X$ be a subcomplex of a metric simplicial complex $X$. 
For every $R\ge 0$, we say that an element
$c\in H_k(X-N_R(K))$ is {\em deep} if it lies in 
$Im(H_k(X-N_{R'}(K))\ra H_k(X-N_R(K)))$ for every $R'\geq R$; equivalently, $c$ 
is deep if belongs to the image of
$$
\underset{\underset{r}{\longleftarrow}}{\lim}\ H_k(X-N_r(K)) \lra H_k(X-N_R(K)). 
$$
We let $H^{Deep}_k(X-N_R(K))$ denote the subgroup of deep homology classes of 
$X-N_R(K)$.  Hence we obtain
an inverse system $\{H^{Deep}_k(X-N_R(K))\}$.  
We say that the deep homology {\em stabilizes} at $R_0$ if the
projection homomorphism 
$$
\underset{\underset{R}{\longleftarrow}}{\lim}\ H^{Deep}_k(X-N_R(K))\ra H^{Deep}_k(X-N_{R_0}(K))$$
is injective.

Specializing the above definition to the case $k=0$,  we arrive at the 
definition of deep  complementary components. If $R\geq 0$, a component $C$ of 
$X-N_R(K)$ is called {\em deep} if it is not contained 
within a finite neighborhood of $K$.  A subcomplex $K$ {\em coarsely separates} $X$ 
if there is an $R$ so that $X-N_R(K)$ has at least two deep components.
A deep component $C$ of $X-N_R(K)$
is said to be {\em stable} if for each $R'\geq R$ the component $C$ 
meets exactly one  deep component of $X- N_{R'}(K)$.  $K$ is said to {\em coarsely
separate} $X$ into (exactly) $m$ components if  there is
an $R$ so that $X-N_R(K)$ consists of exactly $m$ stable deep
components. 

Note that   $H^{Deep}_0(X-N_R(K))$ is freely generated by elements corresponding to 
deep components of $X-N_R(K)$. The deep homology  $H^{Deep}_0(X-N_R(K))$ stabilizes 
at $R_0$ if and only if all deep components of  $X-N_{R_0}(K)$ are stable. 

If $G\acts X$ is a simplicial action of a group on a metric
simplicial complex, then one orbit $G(x)$ coarsely separates $X$ if and only 
if every
$G$-orbit coarsely separates $X$; hence we may simply say that
$G$ coarsely separates $X$.  If $H$ is a subgroup of a finitely
generated group $G$, then we say that $H$ coarsely separates
$G$ if $H$ coarsely separates some (and hence any) Cayley
graph of $G$.

Let $Y, K$ be subcomplexes of a metric simplicial complex $X$. We say that 
$Y$ coarsely separates $K$ in $X$ if there is $R>0$ and two distinct 
 components $C_1, C_2\subset X- N_R(Y)$ so that 
the distance function $d_Y(\cdot)\defeq d(\cdot, Y)$ is unbounded on both 
$K\cap C_1$ and $K\cap C_2$. The subcomplex $Y$ will coarsely separate 
$X$ in this case.

\subsection{Group theoretic preliminaries} 
\label{groupprelim}

{\bf Resolutions, cohomology and relative cohomology.}   
    Let $G$ be  group and $K$ be an Eilenberg-MacLane space for $G$.
If ${\cal M}$ is a system of local coefficients on $K$, then
we have homology and cohomology groups of $K$ with coefficients
in ${\cal M}$:  $H_*(K;{\cal M})$ and $H^*(K;{\cal M})$.
 Now let $A$ be  a $\Z G$-module.  We recall that  a 
{\em resolution of  $A$} is
an exact sequence of $\Z G$-modules:
$$\ldots\to P_n \to \ldots \to P_0 \to A\to 0.$$
 Every $\Z G$-module has a unique projective resolution
up to chain homotopy equivalence.
If $M$ is a $\Z G$-module, then  the 
{\em cohomology of $G$ with coefficients in $M$, $H^*(G;M)$,} is defined
as the homology of chain complex $Hom_{\Z G}(P_*,M)$
where $P_*$ is a projective resolution of the trivial $\Z G$-module
$\Z$; the {\em homology of $G$ with coefficients in $M$, $H_*(G;M)$,}
 is the homology of the chain complex $P_*\otimes_{\Z G}M$.
Using the 1-1 correspondence between $\Z G$-modules $M$ and local
coefficient systems ${\cal M}$ on an Eilenberg-MacLane space $K$, 
we get natural isomorphisms $H_*(K;{\cal M})\simeq H_*(G; M)$
and $H^*(K;{\cal M})\simeq H^*(G; M)$.  Henceforth we will use the same notation
to denote $\Z G$-modules and the corresponding local systems on 
$K(G,1)$'s.

\medskip
\no
{\bf Group pairs.} We now discuss relative (co)homology following \cite{Bieri-Eckmann2}.
   Let $G$ be a group, and 
${\cal H}\defeq \{H_i\}_{i\in I}$
an indexed collection of (not necessarily distinct) subgroups.   We refer to $(G,{\cal H})$
as a {\em group pair}.   Let 
$\amalg_i\, K(H_i,1)\stackrel{f}{\ra} K(G,1)$ be the map 
induced by the inclusions $H_i\ra G$, and
let $K$ be the mapping cylinder of $f$.   We therefore have  a pair
of spaces $(K,\amalg_i \,K(H_i,1))$ since the domain of a  map naturally
embeds in the mapping cylinder.   Given any $\Z G$-module $M$, we define
the relative cohomology $H^*(G,{\cal H};M)$ (respectively  homology
$H_*(G,{\cal H};M)$) to be the cohomology (resp. homology) of the
pair $(K,\amalg_i K(H_i,1))$ with coefficients in the local system $M$. 
As in the absolute case, one can compute relative (co)homology groups
using projective resolutions, see \cite{Bieri-Eckmann2}. For each $i\in I$, let 
$$
\ldots\to Q_n(i) \to \ldots \to Q_0(i) \to \Z\to 0 \quad 
$$
 be  a 
resolution of $\Z$ by projective $\Z H_i$-modules, and let 
$$\ldots\to P_n \to \ldots \to P_0 \to \Z\to 0 \quad $$
be a resolution of $\Z$ by projective $\Z G$-modules.  The inclusions
$H_i\ra G$ induce $\Z H_i$-chain mappings $f_i:Q_*(i)\ra P_*$,  unique up to
chain homotopy.     We define a $\Z G$-chain complex $Q_*$
to be $\oplus_i(\Z G\otimes_{\Z H_i}Q_*(i))$ with an augmentation
$$Q_0\ra  \oplus_i(\Z G\otimes_{\Z H_i}\Z)$$  
induced by the augmentations $Q_0(i)\ra\Z$; the chain mappings $f_i$
yield a $\Z G$-chain mapping $f:Q_*\ra P_*$.    We let $C_*$ be the 
algebraic mapping cylinder of $f$:  this is the chain complex with
$C_i\defeq P_i\oplus Q_{i-1}\oplus Q_i$ with the boundary homomorphism
given by 
\begin{equation}
\label{algcyl}
\D(p_i,q_{i-1},q_i)=(\D p_i +f(q_{i-1}),-\D q_{i-1}  ,\D q_i+q_{i-1}).
\end{equation}
We note that each $C_i$ is clearly projective, a copy $D_*$ of $Q_*$
 naturally sits in $C_*$ as the third summand, and the quotient $C_*/D_*$
is a chain complex of projective $\Z G$-modules.   Proposition 1.2 of
\cite{Bieri-Eckmann2} implies that the relative homology (resp. cohomology) of the 
group pair $(G,{\cal H})$ with coefficients in a $\Z G$-module $M$
(defined as above using local systems on Eilenberg-MacLane spaces) 
is canonically isomorphic to homology of the chain complex
$(C_*/D_*)\otimes_{\Z G}M$  (resp.  $Hom_{\Z G}((C_*/D_*),M)$).

\medskip
\no
{\bf Finiteness properties of groups.} 
The (cohomological) dimension $dim(G)$ of a 
group $G$ is $n$ if $n$ is the minimal integer such that there 
exists a resolution of $\Z$ by  projective $\Z G$-modules: 
$$
0 \to P_n \to ... \to P_0 \to \Z \to 0 .
$$
Recall that $G$ has cohomological dimension $n$ if and only if $n$ is the minimal integer so that
$H^k(G, M)=0$ for all $k>n$ and all $\Z G$-modules $M$. Moreover, if $dim(G)<\infty$ then
$$
dim(G)= \sup \{ n \mid\mbox{$H^n(G, F) \ne 0$ for some free $\Z G$-module $F$}\},
$$
see \cite[Ch. VIII, Proposition 2.3]{Brown}.
If 
$$
1\to G_1 \to G \to G_2 \to 1
$$
is a short exact sequence then $dim(G)\le dim(G_1)+ dim(G_2)$, 
\cite[Ch. VIII, Proposition 2.4]{Brown}. 
If $G'\subset G$ is a subgroup then $dim(G')\le dim(G)$. 

A {\em partial resolution} of a $\Z G$-module $A$ is an exact sequence
 $\Z G$-modules:
$$ P_n \to \ldots \to P_0 \to A\to 0 .$$
If $A_*$: 
$$ ... \to A_n \to A_{n-1}\to \ldots \to A_0 \to A\to 0 $$
is a chain complex then we let $[A_*]_n$ denote the {\em $n$-truncation 
of $A_*$}, i.e. 
$$
A_n \to \ldots \to A_0 \to A\to 0 .  
$$
A group $G$ is of type $FP_n$ if there exists a partial 
resolution of $\Z$ by finitely generated projective $\Z G$-modules: 
$$
P_n \to ... \to P_0 \to \Z \to 0 .
$$
The group $G$ is of type $FP$ (resp. $FL$) if there exists 
a finite resolution of 
$\Z$ by finitely generated projective (resp. free) $\Z G$-modules. 
We will also refer to groups of type $FP$ as {\em groups of finite type}.

\begin{lemma}
\label{shortres}
1. If $G$ is of type $FP$ then $dim(G)=n$ if and only if
$$
n= \max \{ i : H^i(G,\Z G) \ne 0\}.
$$
2. If $dim(G)=n$ and $G$ is of type $FP_n$ then there exists a 
resolution of $\Z$ by finitely generated projective $\Z G$-modules: 
$$
0\to P_n \to ... \to P_0 \to \Z \to 0  .
$$
In particular $G$ is of type $FP$.
\end{lemma}

\proof The first assertion follows from \cite[Ch. VIII, Proposition 5.2]{Brown}. 
We prove 2. Start with a partial resolution
$$
P_n \to P_{n-1} \to ... \to P_0 \to \Z \to 0
$$
where each $P_i$ is finitely generated projective. 
By \cite[Ch. VIII, Lemma 2.1]{Brown},
the kernel $Q_n:= \ker[P_{n-1} \to P_{n-2}]$ is projective. 
However $P_n$ maps onto $Q_n$, hence $Q_n$ is also finitely generated. Thus 
replacing $P_n$ with $Q_n$ we get the required resolution. \qed

Examples of groups of type $FP$ and $FL$ are given by fundamental groups of finite 
Eilenberg-MacLane complexes, or more generally, groups acting freely cocompactly 
on acyclic complexes. According to the Eilenberg-Ganea theorem, if $G$ is a finitely 
presentable group of type $FL$ then $G$ admits a finite $K(G,1)$ of 
dimension $\max(dim(G),3)$. 

Let $G$ be a group, let ${\cal H}\defeq\{H_i\}_{i\in I}$ be an indexed  collection
of subgroups, and  let $$\eps:\oplus_i \,(\Z G\otimes_{\Z H_i}\Z)\ra \Z$$
be induced by the usual augmentation $\Z G\ra \Z$.     Then the 
group pair {\em $(G,{\cal H})$ has finite
type} if the $\Z G$-module $Ker(\eps)$ admits
a finite length resolution by finitely generated projective $\Z G$-modules.
If the index set $I$ is finite and the groups $G$ and $H_i$ are of type $FP$ then
the group pair $(G,{\cal H})$ is of finite type, and one obtains the
desired resolution of $Ker(\eps)$ using the quotient $C_*/D_*$ where
$(C_*,D_*)$ is the pair given by the algebraic mapping cylinder
construction (\ref{algcyl}).

\medskip
For the next three topics, the reader may consult 
\cite{bieri,Bieri-Eckmann1,Bieri-Eckmann2,Brown,dunwoodydicks}.

\medskip
\no
{\bf Duality groups. }
Let $G$ be a group of type $FP$.  Then $G$ is an {\em $n$-dimensional duality group}
if $H^i(G;\Z G)=\{0\}$ when $i\neq n=dim(G)$, and $H^i(G;\Z G)$ is 
torsion-free, \cite{Bieri-Eckmann1}.  There is an alternate definition of duality groups 
involving isomorphisms $H^i(G;M)\simeq H_{n-i}(G;D\otimes M)$ for 
a suitable dualizing module $D$ and arbitrary $\Z G$-modules
$M$, see \cite{Bieri-Eckmann1,Brown}. Examples of duality
groups include:

1. The fundamental groups of compact aspherical manifolds with aspherical
boundary, where the inclusion of each boundary component induces
a monomorphism of fundamental groups.

2. Torsion-free S-arithmetic groups, 
\cite{Bieri-Eckmann1,BS1}.

3.  $2$-dimensional one-ended groups of type $FP_2$ 
\cite[Proposition 9.17]{bieri}; for instance
torsion-free, one-ended, one-relator groups.  

4.  Any group which can act freely, cocompactly, and simplicially
on an acyclic simplicial complex $X$, where $H^i_c(X)$ vanishes
except in dimension $n$, and $H^n_c(X)$ is torsion-free.

\medskip
\no
{\bf Poincar\'e duality groups.}  These form a special class of duality groups.
If $G$ is an $n$-dimensional duality group and $H^n(G;\Z G)=\Z$, then
$G$ is an {\em $n$-dimensional Poincare duality group} ($PD(n)$ group).   
As in the case of duality groups, there is an alternate definition
involving isomorphisms $H^i(G;M)\simeq H_{n-i}(G;D\otimes M)$ where 
$M$ is an arbitrary $\Z G$-module and the orientation $\Z G$-module $D$ is 
isomorphic to $\Z$  as an abelian group. Examples include:

1. Fundamental groups of closed aspherical manifolds.

2. Fundamental groups of aspherical finite Poincare complexes.  Recall that
an (orientable)  {\em Poincare complex of formal dimension $n$} is a finitely 
dominated complex $K$
together with a fundamental class $[K]\in H_n(K;\Z)$ so that the
cap product operation $[K]\cap: H^k(K;M)\ra H_{n-k}(K;M)$
is an isomorphism for every local system $M$ on $K$ and for $k=0,\ldots,n$.

3.  Any group which can act freely, cocompactly, and simplicially
on an acyclic simplicial complex $X$, where $X$ has the same
compactly supported cohomology as $\R^n$.

4. Each torsion-free Gromov-hyperbolic group $G$ whose boundary is a homology 
manifold with the homology of sphere (over $\Z$),  see \cite{bm}. Note 
that every such group is the  fundamental group of a  finite  aspherical 
Poincare complex, namely the $G$-quotient of a Rips complex of $G$.    

Below are several useful facts about Poincare duality groups 
(see \cite{Brown}):

(a) If $G$ is a $PD(n)$ group and $G'\subset G$ is a subgroup 
then $G'$ is a $PD(n)$ group if and only if the index $[G:G']$ is finite. 

(b) If $G$ is a $PD(n)$ group which is contained in a torsion-free 
group $G'$ as a finite index subgroup, 
then $G'$ a $PD(n)$ group.

(c) If $G\times H$ is a $PD(m)$ group then 
$G$ and $H$ are $PD(n)$ and  $PD(k)$ groups, where $m=n+k$. 

(d) If $G\rtimes H$ is a semi-direct product where $G$
is a $PD(n)$-group and $H$ is a $PD(k)$-group, then $G\rtimes H$
 is a $PD(n+k)$-group. 
See \cite[Theorem 3.5]{Bieri-Eckmann1}. 

There are several questions about $PD(n)$ groups and their 
relation with fundamental groups of aspherical manifolds.
It was an open question going back
to Wall \cite{Wall} whether every $PD(n)$ group is the fundamental group
of a closed aspherical manifold.
The answer to this is yes in dimensions $1$ and $2$,
\cite{stallings,EL,EM}.    Recently, Davis in \cite{davis}
gave examples for $n\geq 4$ of  $PD(n)$ groups which do not admit a finite
presention, and these groups are clearly not fundamental groups
of compact manifolds.  This leaves open several questions: 

1.  Is every finitely
presented $PD(n)$ group the fundamental group of a compact
aspherical manifold?

2.  A weaker version of 1: Is every finitely
presented $PD(n)$ group the fundamental group of a finite
aspherical complex?  Equivalently, by Eilenberg-Ganea, one
may ask if every such group is of type $FL$.

3. Does every $PD(n)$ group act freely and cocompactly
on an acyclic  complex?  We believe
this question is  open for groups of type $FP$.  One
can also ask if every $PD(n)$ group acts freely and cocompactly
on an acyclic $n$-manifold.

\medskip
\no
{\bf Poincare duality pairs.}
Let $G$ be an $(n-1)$-dimensional group of type $FP$, 
and let $H_1,\ldots,H_k\subset G$ be $PD(n-1)$ subgroups
of $G$.  Then the group pair $(G;H_1,\ldots,H_k)$ is an {\em $n$-dimensional
Poincare duality pair, or  $PD(n)$ pair},
if the double of $G$ over the $H_i$'s is a $PD(n)$ group.
We recall that the double of $G$ over the $H_i$'s is the
fundamental group of the graph of groups ${\cal G}$,
where ${\cal G}$ has two vertices labelled by $G$,
$k$ edges with the $i^{th}$  edge labelled by  $H_i$, and edge
monomorphisms are the inclusions $H_i\ra G$.
An alternate homological 
definition of $PD(n)$ pairs is the following: a group pair
$(G,\{H_i\}_{i\in I})$ is a $PD(n)$ pair if it has finite
type, and $H^*(G,\{H_i\};\Z G)\simeq H^*_c(\R^n)$. 
For a discussion of these and other equivalent definitions, 
see \cite{Bieri-Eckmann2,dunwoodydicks}.  We will sometimes
refer to the system of subgroups $\{H_i\}$ as the 
{\em peripheral structure} of the $PD(n)$ pair, and the 
$H_i$'s as peripheral subgroups.
The first class of examples of duality groups mentioned
above have natural peripheral structure which makes
them $PD(n)$ pairs. In \cite{KK}  we proved that if $G$ 
is a torsion-free Gromov-hyperbolic group whose boundary is 
homeomorphic to the  Sierpinski carpet $S$,  then $(G; H_1,...,H_k)$ 
is  a $PD(3)$ group pair, where $H_i$'s are representatives of 
conjugacy classes of stabilizers of the peripheral circles of $S$ in $\geo G$.   
If $(G;H_1,\ldots,H_k)$ 
is a  $PD(n)$ pair, where $G$ and each 
$H_i$ admit  a finite Eilenberg-MacLane space $X$ and $Y_i$
respectively, then the inclusions $H_i\ra G$ induce
a map $\sqcup_i Y_i\ra X$ (well-defined up to homotopy)
whose mapping cylinder $C$ gives a {\em Poincare pair}
$(C;\sqcup_i Y_i)$, i.e. a pair which satisfies Poincare duality 
for manifolds with boundary with local coefficients 
(where $\sqcup_i Y_i$ serves as the boundary of $C$). Conversely, if $(X,Y)$ is a Poincare pair
where $X$ is aspherical and  
$Y$ is a union of aspherical components $Y_i$, 
then $(\pi_1(X);\pi_1(Y_1),\ldots,\pi_1(Y_k))$
is a  $PD(n)$ pair.

\begin{lemma}
\label{hi'sdistinct}
Let $(G,\{H_i\})$ be a $PD(n)$ pair, where $G$ is not a 
$PD(n-1)$ group.  Then the subgroups $H_i$ are pairwise
non-conjugate maximal $PD(n-1)$ subgroups.
\end{lemma}
\proof
If $H_i$ is conjugate to $H_j$ for some $i\neq j$, then the double $\hat G$ of $G$ over the 
peripheral subgroups would contain an infinite index subgroup isomorphic 
to the $PD(n)$ group $H_i\times \Z$.  The group $\hat G$ is a $PD(n)$ group, which 
contradicts property (a) of Poincare duality groups listed above.

We now prove that each $H_i$ is maximal. 
 Suppose that $H_i\subset H\subset G$, where $H\ne H_i$ is a 
$PD(n-1)$ group. Then $[H: H_i]< \infty$. Pick $h\in H- H_i$. Then there 
exists a finite index subgroup $F_i\subset H_i$ which is normalized 
by $h$. Consider the double $\hat{G}$ of $G$ along the collection 
of subgroups $\{H_i\}$, and let $\hat{G}\acts T$ be the associated action on 
the Bass-Serre tree. Since $G$ is not a $PD(n-1)$ group, $H_i\neq G$
for each $i$, and so there is a unique vertex $v\in T$ fixed by $G$.
The involution of the graph of groups defining $\hat G$ induces
an involution of $\hat G$ which is unique up to
an inner automorphism;  let $\tau:\hat G\ra\hat G$  be an induced
involution which fixes $H_i$ elementwise.  Then $G'\defeq \tau(G)$
fixes a vertex $v'$ adjacent to $v$, where the edge 
$\ol{vv'}$ is fixed by $H_i$.   So $h'\defeq\tau(h)$ belongs to $\tau(G)=G'$
but $h'$ does not fix $\ol{vv'}$.  Therefore the fixed point sets of
$h$ and $h'$ are disjoint, which implies that $g\defeq hh'$ acts
on $T$ as a hyperbolic automorphism.  Since $h'\in Normalizer(\tau(F_i))
=Normalizer(F_i)$, we get $g\in Normalizer(F_i)$.  Hence
the subgroup $F$ generated by $F_i$ and $g$  is a semi-direct
product $F=F_i\rtimes \<g\>$, and $\<g\>\simeq\Z$ since
$g$ is hyperbolic.   The group $F$ is a $PD(n)$ group (by property (d))  
sitting as an infinite index subgroup of the $PD(n)$ group $G$,
which contradicts property (a).
\qed

\subsection{Algebraic preliminaries}
\label{algprelim}

In this section we introduce a notion of ``morphism'' between inverse systems.
Approximate isomorphisms, which figure prominently in the remainder of the
paper, are maps between inverse (or direct) systems which fail to be isomorphisms in
a controlled way, and for many purposes are as easy to work with as 
isomorphisms.

\medskip
\no
{\bf Approximate morphisms between inverse and direct systems.} 
Recall that a partially ordered set $I$ is  {\em directed} if for each $i, j\in I$ there exists 
$k\in I$ such that $k\ge i, j$.  An inverse system of (abelian) groups indexed 
by a directed set $I$ is a collection of abelian groups $\{A_i\}_{ i\in I}$ and 
homomorphisms ({\em projections}) $p^{j}_{i}:A_{i}\ra A_{j}$, $i\ge j$ so that 
$$
\mbox{$p^i_i=id$ and  $p^k_j \circ p^j_i = p^k_i$}$$
 for any $i\le j\le k$. 
(One may weaken these assumptions but they will suffice for our purposes.) 
 We will often denote the inverse system by $(A_\bul,p_\bul)$ or $\{ A_i\}_{i\in I}$.  
 Recall that a subset $I'\subset I$ of a partially 
ordered set is  {\em cofinal} if  for every $i\in I$ there is an 
$i'\in I'$ so that $i'\geq i$. 

Let $\{A_i\}_{i\in I}$ and $\{B_j\}_{j\in J}$ be two inverse
systems of (abelian) groups indexed by $I$ and $J$, with the projection maps
$p^{i'}_{i}:A_{i}\ra A_{i'}$ and $q^{j'}_{j}:B_{j}\ra B_{j'}$. 
The directed sets appearing later
in the paper will be order isomorphic to $\Z_+$ with the usual order. 

\begin{definition}
Let $\al$ be an order preserving, partially defined, map from $I$ to $J$.
Then $\al$ is {\em cofinal} if it is defined on a subset of the 
form $\{i\in I\mid i\geq i_0\}$ for some $i_0\in I$, and 
 the image of every
cofinal subset $I'\subset I$ 
is a cofinal subset $\al(I')\subset J$. 
\end{definition}

\begin{definition}
\label{approxmorphdef}
Let $\al:I\ra J$ be a  cofinal map.  Suppose  that    
$(\{A_i\}_{i\in I},p_\bul)$ and $(\{B_j\}_{j\in J},q_\bul)$   are  
inverse systems.  
Then a family of homomorphisms $f_i:A_i\ra B_{\al(i)}$, $i\in I$,
is an {\em $\al$-morphism from $\{A_i\}_{i\in I}$
to $\{B_j\}_{j\in J}$} if 
\begin{equation}
 \label{f'scompat}
q_{\al(i)}^{\al(i')}\circ f_i=f_i\circ p_i^{i'}
\end{equation}
whenever $i,i'\in I$ and $i\geq i'$. The {\em saturation} 
$\hat{f}_\bul^\bul$ of the $\al$-morphism $f_\bul$ 
is the collection of maps $\hat f_i^j: A_i \to B_j$ of the form
$$
q_{\al(k)}^j\circ f_k \circ p_i^k\  .
$$ 
In view of (\ref{f'scompat}) this definition is consistent, and 
$\hat f_\bul^\bul$ is compatible with the projection maps
of $A_\bul$ and $B_\bul$. 
\end{definition}

\no  Suppose that $\{A_i\}_{i\in I}$, $\{B_j\}_{j\in J}$, $\{C_k\}_{k\in K}$ 
are inverse systems, $\al:  I\to J$, $\be: J\to K$ are cofinal  
maps. Then the composition of $\al$- and $\be$-morphisms 
$$
f_\bul : A_\bul \to B_\bul, \quad g_\bul : B_\bul \to C_\bul
$$
is a $\ga$-morphism for the cofinal map $\ga= \be \circ \al :  I\to K$. 
(The composition $\be\circ\al$ is defined on the subset 
$Domain(\al)\cap \al^{-1}(Domain(\beta))$
which contains 
$\{i : i\ge i_1\}$ where $i_1$ is an upper bound for non-cofinal subset 
$\al^{-1}(J-Domain(\beta))$ in $I$.)

\begin{definition} Let $A_\bul \stackrel{f_\bul}{\ra} B_\bul$ 
be an $\al$-morphism of inverse systems $(A_\bul , p_\bul), (B_\bul , q_\bul)$.

1. When $I$ is totally ordered, we define  $Im(\hat{f}_\bul^j)$,  the  {\em  image}  of  $f_\bul$  
in  $B_j$, to be  $\cup\{ Im(\hat f_i^j:A_i\ra B_j)\mid \al(i)\geq j\}$.  

2.  Let $\om:I\ra I$ be a function. Then $f_\bullet$ is an {\em $\om$-approximate monomorphism} 
if for every $i\in I$ we have
$$Ker(A_{\om(i)}\stackrel{f_{\om(i)}}{\lra} B_{\al(\om(i))})
\subset Ker(A_{\om(i)}\stackrel{p_\bul}{\lra}A_i).$$

3. Suppose $I$ is totally ordered.  If  $\bar\om:J\ra J$ is a function,
then  $f_\bullet$ is an {\em $\bar\om$-approximate epimorphism} 
if for every $j\in J$ we have:
$$
Im(B_{\bar\om(j)}\stackrel{q_\bul}{\lra}B_j)\subset Im(\hat{f}_\bul^j).
$$

4. Suppose $I$ is totally ordered.
If $\om:I\ra I$ and $\bar\om:J\ra J$ are functions,
then $f$ is an {\em $(\om,\bar\om)$-approximate isomorphism} 
if  both 2 and 3 hold.
\end{definition}

We will frequently suppress the functions $\al,\,\om,\,\bar\om$  when
speaking of morphisms, approximate monomorphisms (epimorphisms,
isomorphisms). Note that the inverse limit of an approximate 
monomorphism  (epimorphism, isomorphism)  is a  monomorphism  (epimorphism,
isomorphism) of inverse limits.

Note that an $\al$-morphism induces a homomorphism between inverse limits, 
since for each cofinal subset $J'\subset J$ we have:
$$
\underset{\underset{j\in J}{\longleftarrow}}{\lim} \ B_j 
\cong \underset{\underset{j\in J'}{\longleftarrow}}{\lim}\ B_j\ .
$$
Similarly, an approximate isomorphism of inverse systems induces an 
isomorphism of their inverse limits. However the converse is not true. 
For instance, let $A_i:= \Z$ for each $i\in \N$, where $\N$ has the 
usual order. Let 
$$
p^{i-n}_i: A_i\to A_{i-n} \hbox{~~be the index $n$ inclusion}.  
$$ 
It is clear that the inverse limit of this system is zero. We leave it to
the reader to verify that the system $(A_\bul,p_\bul)$ is not approximately 
isomorphic to zero inverse system. 


We have similar definitions for homomorphisms of direct systems. 
A direct system of (abelian) groups indexed 
by a directed set $I$ is a collection of abelian groups $\{A_i\}_{ i\in I}$ and 
homomorphisms ({\em projections}) $p_i^j:A_{i}\ra A_{j}$, $i\le j$ so that 
$$
p^i_i=id, \quad p_j^k\circ p_i^j=p_i^k$$
 for any $i\leq j\leq k$. 
We often denote the direct system by $(A_\bul,p_\bul)$.  
Let $\{A_i\}_{i\in I}$ and $\{B_j\}_{j\in J}$ be two direct 
systems of (abelian) groups indexed by directed sets 
$I$ and $J$, with projection maps
$p_i^{i'}:A_{i}\ra A_{i'}$ and $q_{j}^{j'}:B_{j}\ra B_{j'}$. 

\begin{definition}
Let $\al:I\ra J$ be a  cofinal map.
Then a family of homomorphisms $f_i:A_i\ra B_{\al(i)}$, $i\in I$,
is a {\em $\al$-morphism of the direct systems $\{A_i\}_{i\in I}$
and $\{B_j\}_{j\in J}$} if 
$$
q_{\al(i)}^{\al(i')}\circ f_i=f_{i'}\circ p_i^{i'}
$$
whenever $i\leq i'$.
We define the saturation $\hat{f}_\bul^\bul$ the same way as
 for morphisms of inverse systems. 
\end{definition}

\begin{definition} 
Let $f_\bullet:A_\bullet \ra B_\bullet$ be an $\al$-morphism of direct systems: 
$$
f_\bullet=\{f_i : A_i \to B_{\al(i)}, i\in I\}.$$

1.  When $I$ is totally ordered
we define  $Im(\hat{f}_\bul^j)$,  the  {\em  image}  of  $f_\bul$  
in  $B_j$, to be  $\cup\{Im(\hat f_i^j)\mid \al(i)\leq j\}$.

2.  Let $\om:I\ra I$ be a function. Then $f_\bullet$ is an {\em $\om$-approximate monomorphism} 
if for every $i\in I$ we have 
$$
Ker(A_{i}\stackrel{f_{i}}{\lra} B_{\al(i)})
\subset Ker(A_{i}\stackrel{p_\bul}{\lra}A_{\om(i)}).$$

3. Suppose $I$ is totally ordered, and  $\bar\om:J\ra J$ is a function.  
$f_\bullet$ is an {\em $\bar\om$-approximate epimorphism} 
if for every $j\in J$ we have:
$$
Im(B_{j}\stackrel{q_{\bul}}{\lra}B_{\bar\om(j)})\subset
Im(\hat{f}_\bul^{\bar\om (j)}).$$

4. Suppose $I$ is totally ordered and  $\om:I\ra I$ and $\bar\om:J\ra J$
are functions.
Then $f$ is an {\em $(\om,\bar\om)$-approximate isomorphism} 
if  both 2 and 3 hold.
\end{definition}

 An inverse (direct) system $A_\bullet$ is 
said to be constant if $A_i=A_j$ and $p^i_j=id$ for each $i, j$. An 
inverse (direct) system $A_\bullet$ is {\em approximately constant} if 
there is an  approximate isomorphism between it and a constant system
(in either direction).   Likewise, an inverse or direct system is {\em
approximately zero} if it is approximately isomorphic to a zero
system.

The proof of the following lemma is straightforward and is left to the reader. 

\begin{lemma}
The composition of two approximate monomorphisms (epimorphisms,
isomorphisms) is an approximate monomorphism  (epimorphism, isomorphism).
\end{lemma}

\bigskip
The remaining material in this section relates to the category theoretic
behavior of approximate morphisms and a comparison with pro-morphisms,
and it will not be used elsewhere in the paper.

In what  follows  $(A_\bul,p_\bul)$ and $(B_\bul,q_\bul)$ 
will once again denote inverse systems indexed
by $I$ and $J$ respectively.   However, for simplicity we will assume that
$I$ and $J$ are both totally ordered.

\begin{definition}
Let $f_\bul:A_\bul\ra B_\bul$ be an $\al$-morphism with  
saturation $\hat f_\bul^\bul$.
The {\em kernel} of $f_\bul$ is the inverse system $\{K_i\}_{i\in I}$ where
$K_i\defeq Ker(f_i:A_i\ra B_{\al(i)})$ with the projection maps obtained from
the projections of $A_\bul$ by restriction.  We define the {\em image} of 
$f_\bul$ to be the inverse system $\{D_j\}_{j\in J}$
where $D_j\defeq Im(\hat f^j_\bul)$,
with the projections coming from the projections of $B_\bul$. Note that $D_j$  
is a  subgroup  of  $B_j$, $j\in J$.  We
also define the {\em cokernel $coKer(f_\bul)$ of} $f_\bul$,    
as the inverse  system $\{C_j\}_{j\in J}$ where $C_j\defeq B_j/D_j$.
\end{definition}

An inverse (respectively direct) system of abelian groups $A_\bul$ is {\em pro-zero}  
if for every $i\in I$ there exists $j\ge i$ such that $p_j^i : A_j \to A_i$ 
(resp. $p_i^j : A_i \to A_j$) is  zero  (see \cite[Appendix 3]{Artin-Mazur}). 
Using this language we may reformulate  the definitions of 
approximate monomorphisms:

\begin{lemma}
Let $f_\bul : A_\bul \to B_\bul$ be a morphism of inverse systems of 
abelian groups. Then 

1. $f_\bul$ is an approximate monomorphism iff its kernel 
$K_\bul := Ker(f_\bul)$ is pro-zero. 

2.   $f_\bul$ is an approximate 
epimorphism iff its cokernel is a pro-zero inverse system. 

3.  $f_\bul$ is an approximate isomorphism iff both
$Ker(f_\bul)$ and $coKer(f_\bul)$ are pro-zero systems.
\end{lemma}
\proof This is immediate from the definitions.\qed

For a fixed cofinal map $\al:I\ra J$, the collection of $\al$-morphisms from
 $A_\bul$ to $B_\bul$ forms an abelian group the obvious way.  In order to
compare morphisms $A_\bul\ra B_\bul$ with different index maps
$I\ra J$, we introduce an equivalence relation:

\begin{definition}
\label{approxequivdef}
Let $f:A_\bul\ra B_\bul$ and $g:A_\bul\ra B_\bul$ be morphisms
with saturations $\hat  f_\bul^\bul$ and $\hat g_\bul^\bul$.
Then {\em $f_\bul$ is  equivalent $g_\bul$} if there is a cofinal function
$\rho:J\ra I$ so that for all $j\in J$, both 
$\hat f_{\rho(j)}^j$ and $\hat g_{\rho(j)}^j$ are defined, and they coincide.
\end{definition}
\no
This equivalence relation is compatible with composition of approximate
morphisms.   Hence we obtain a category
$Approx$ where the objects are inverse systems of abelian groups and  the morphisms
are equivalence classes of approximate morphisms.   An {\em approximate
inverse} for an approximate morphism $f_\bul$ is an approximate
morphism $g_\bul$ which inverts $f_\bul$ in $Approx$.  

\begin{lemma}
\label{bitsandpieces}
  Suppose $I,J\cong\Z_+$,
$D_\bul$ is a sub inverse system of $A_\bul$  (i.e.  $D_i\subset  A_i$,  
$i\in  I$), and let
$Q_\bul$ be the quotient system: $Q_i\defeq A_i/D_i$.  Then 

1. The
morphism $A_\bul\ra Q_\bul$ induced by the canonical epimorphisms
$A_i\ra Q_i$ has an approximate inverse iff $D_\bul$ is a pro-zero system.

2. The morphism $D_\bul\ra A_\bul$ defined by the inclusion homomorphisms
$D_i\ra A_i$ has an approximate inverse iff $Q_\bul$ is a pro-zero system.

3.  If $f_\bul:A_\bul\ra B_\bul$ is a morphism, $Ker(f_\bul)$ is zero (i.e. $Ker(f_\bul)_i=\{0\}$
for all $i\in I$), and $Im(f_\bul)=B_\bul$, then $f_\bul$ has an approximate inverse.
\end{lemma}
\proof  We leave the ``only if'' parts of 1 and 2 to the reader.

When $D_\bul$ is pro-zero  the  map   $\be:I\ra I$ defined by 
$$
\be(i)\defeq \max\{i'\mid D_i\subset Ker(A_i\ra A_{i'})\}$$
is cofinal.  Let 
$g_\bul:Q_\bul\ra A_\bul$ be the $\be$-morphism where
$g_i:A_i/D_i=Q_i\ra A_{\be(i)}$ is induced by the projection
$A_i\ra A_{\be(i)}$.  One checks that $g_\bul$ is an approximate
inverse  for $A_\bul\ra Q_\bul$.

Suppose $Q_\bul$ is pro-zero.   Define a cofinal map $\be:I\ra I$
by 
$$
\be(i)\defeq \max\{i'\mid Im(A_i\ra A_{i'})\subset D_{i'}\},$$ 
 and
let $g_\bul: A_\bul\ra D_\bul$ be the $\be$-morphism where
$g_i:A_i\ra D_{\be(i)}$ is induced by the projection $A_i\ra A_{\be(i)}$.
Then $g_\bul$ is an approximate inverse for the inclusion $D_\bul\ra A_\bul$.

Now suppose $f_\bul:A_\bul\ra B_\bul$ is an $\al$-morphism with zero
kernel and cokernel.  Let $J'\defeq \al(I)\subset J$, and define
$\be':J'\ra I$ by $\be'(j)=\min \al^{-1}(j)$.   Define  a  cofinal 
map $\si:J\ra J'$
by $\si(j)\defeq \max\{j'\in J'\mid j'\leq j\}$;  let $\be:J\ra I$
be the composition $\be'\circ\si$, and define a $\be$-morphism
$g_\bul$ by $g_j\defeq f^{-1}_{\be(j)}\circ q_j^{\si(j)}$.  Then
$g_\bul$ is the desired approximate inverse for $f_\bul$.
\qed

\begin{lemma}
\label{isoiffapproxinvertible}
 Let $f_\bul:A_\bul\ra B_\bul$ be a  morphism.  

1. If $f_\bul$ has an approximate inverse then it is an 
approximate isomorphism.

2. If $f_\bul$ is an approximate isomorphism and $I,J\cong\Z_+$ then
$f_\bul$ has an approximate inverse.
\end{lemma}
\proof
Let $f_\bul:A_\bul\ra B_\bul$ and $g_\bul:B_\bul\ra A_\bul$ be 
$\al$ and $\be$ morphisms respectively, and let $g_\bul$
be an approximate inverse for $f_\bul$.  Since $h_\bul\defeq g_\bul\circ f_\bul$
is equivalent to $id_{A_\bul}$ then for all $i$ there is an $i'\geq i$
so that $\hat h_{i'}^i$ is defined and $\hat h_{i'}^i=p_{i'}^i$.
Letting $\ga\defeq \be\circ\al$ we have, by the definition of the
saturation $\hat h_\bul^\bul$, $p_{i'}^i=\hat h_{i'}^i=p_{\ga(i)}^i\circ
h_{i'}$.  So $Ker(h_{i'})\subset Ker(p_{i'}^i)$.  Thus $f_\bul$ is 
an approximate monomorphism.  The proof that $f_\bul$ is
an approximate epimorphism is similar. 

We now prove part 2.
Let $\{K_i\}_{i\in I}$ be the kernel of $f_\bul$, let 
$\{Q_i\}_{i\in I}=\{A_i/K_i\}_{i\in I}$ be the quotient system,
and let $\{D_j\}_{j\in J}$ be the image of $f_\bul$.  Then
$f_\bul$ may be factored as $f_\bul=t_\bul\circ s_\bul\circ r_\bul$ where
$r_\bul:A_\bul\ra Q_\bul$ is induced by the epimorphisms $A_i\ra A_i/K_i$,
$s_\bul:Q_\bul\ra D_\bul$ is induced by the homomorphisms of quotients,
and $t_\bul:D_\bul\ra B_\bul$ is the inclusion.   By Lemma \ref{bitsandpieces}, 
$s_\bul$ has an approximate inverse.   When the kernel and cokernel
of $f_\bul$ are pro-zero then $r_\bul$ and $t_\bul$ also admit approximate
inverses by Lemma \ref{bitsandpieces}. 
 Hence $f_\bul$ has an approximate inverse in this case.
\qed

\bigskip
\no
{\bf Relation with Grothendieck's pro-morphisms.}
Below we relate the notions of $\al$-morphisms, approximate monomorphisms 
(epimorphisms, isomorphisms) with Grothendieck's pro-morphisms. Strictly speaking 
this is unnecessary for the purposes of this paper, however it puts our definitions 
into perspective.  Also, readers who prefer the language of pro-categories
may use Lemma \ref{faithful} and Corollary \ref{approxiffpro}
 to translate the theorems of sections 
\ref{cpd} and \ref{cad} into pro-theorems.

\begin{definition}
Let $\{A_i\}_{i\in I},\, \{B_j\}_{j\in J}$ be 
inverse systems.  The group of pro-mor\-phisms $proHom(A_\bul, B_\bul)$ is defined as 
$$
\underset{\underset{j\in J}{\longleftarrow}}{\lim}~ 
 \underset{\underset{i\in I}{\longrightarrow}}{\lim}~ Hom(A_i, B_j)
$$
(see \cite{Grot}, \cite[Appendix 2]{Artin-Mazur}, \cite[Ch II, \S 1]{DS}). 
The {\em identity pro-morphism} is the element of
$proHom(A_\bul,A_\bul)$ determined by 
$(id_{A_j})_{j\in I}\in 
\prod_{j}\,\underset{\underset{i\in I}{\longrightarrow}}{\lim}~ Hom(A_i, A_j)$. 
\end{definition}
\no
This yields a category\footnote{By relaxing the definition of inverse systems, 
this category becomes an abelian category, \cite[Appendix 4]{Artin-Mazur}.  However we will not discuss
this further.} 
{\it Pro-Abelian} where the objects are inverses systems of 
abelian groups and the morphisms are the pro-morphisms. 
 A {\em pro-isomorphism} is an isomorphism in this category.

By the definitions of direct and inverse limits, an element of 
$proHom(A_\bul, B_\bul)$ can be represented by an admissible ``sequence'' 
$$
([h_{\rho(j)}^j : A_{\rho(j)} \to B_j])_{j\in J}
$$
of equivalence classes of homomorphisms $h_{\rho(j)}^j : A_{\rho(j)} \to B_j$; 
here two homomorphisms 
$h_i^j : A_i \to B_j, h_k^j : A_k \to B_j$ are equivalent if there exists $\ell \ge i,  k$ 
such that
$$
h_i^j \circ p_\ell^i =  h_k^j \circ p_\ell^k;
$$
and the ``sequence'' is {\em admissible} if for each 
$j\geq j'$ there is an $i\geq  \max\{\rho(j),\rho(j')\}$ so that
$$
q_j^{j'}\circ h_{\rho(j)}^j\circ p_i^{\rho(j)}=h_{\rho(j')}^{j'}\circ p_i^{\rho(j')}.
$$

Given a cofinal map $\al:I\ra J$ between directed sets, 
we may construct\footnote{Using the axiom of choice 
we  pick $\rho(j)\in  \al^{-1}(j)$.} a function $\rho:J\ra I$
so that $\al(\rho(j))\geq j$ for all $j$; then any $\al$-morphism 
$f_\bul:A_\bul\ra B_\bul$ induces an admissible sequence 
$([\hat f_{\rho(j)}^j:A_{\rho(j)}\ra B_j]\}_{j\in J}$.  The
corresponding element $pro(f_\bul)\in proHom(A_\bul,B_\bul)$ is independent
of the choice of $\rho$  by condition (\ref{f'scompat}) of Definition
\ref{approxmorphdef}.   

\begin{lemma}
\label{faithful}

1.  If $f:A_\bul\ra B_\bul$ and $g:A_\bul\ra B_\bul$  are 
morphisms, then $pro(f)=pro(g)$ iff  $f_\bul$ is equivalent to $g_\bul$.
In other words, $pro$ descends to a faithful functor from 
$Approx$ to {\it Pro-Abelian}.

2.  When $I,\,J\cong\Z_+$ then every pro-morphism from $A_\bul$ to $B_\bul$ arises as
$pro(f_\bul)$ for some approximate morphism $f_\bul:A_\bul\ra B_\bul$. 
Thus $pro$ descends to a fully faithful functor from 
$Approx$ to {\it Pro-Abelian} in this case.
\end{lemma}
\proof
The first assertion follows readily from the definition of $proHom(A_\bul,B_\bul)$
and Definition \ref{approxequivdef}.

Suppose $I,\,J\cong\Z_+$ and  $\phi\in proHom(A_\bul, B_\bul)$  is   
represented by an admissible sequence
$$
([h_{\rho_0(j)}^j : A_{\rho_0(j)} \to B_j])_{j\in J}.
$$
We define $\rho:J\to  I$  and  another admissible
sequence $(\bar h_{\rho(j)}^j:A_{\rho(j)}\ra B_j)_{j\in J}$
representing $\phi$ by setting 
$\rho(0)=\rho_0(0)$, $\bar h_{\rho(0)}^0\defeq h_{\rho_0(0)}^0$, and 
inductively choosing $\rho(j)$, $\bar h_{\rho(j)}^j$ so that
$\rho(j)>\rho(j-1)$, $\bar h_{\rho(j)}^j\defeq h_{\rho_0(j)}^j\circ p_{\rho(j)}^{\rho_0(j)}$
and $q_j^{j-1}\circ \bar h_{\rho(j)}^j=\bar h_{\rho(j-1)}^{j-1}\circ p_{\rho(j)}^{\rho(j-1)}$.   Note  that  the  mapping  $\rho$ is  strictly  
increasing  and  hence  cofinal.  
Now define a cofinal map $\al:\Z_+\ra\Z_+$ by setting
$\al(i)\defeq\max\{j\mid \rho(j)\leq i\}$ for $i\geq\rho(0)=\rho_0(0)$.
We then  get an $\al$-morphism $f_\bul:A_\bul\ra B_\bul$ where
$f_i\defeq\bar h_{\rho(\al(i))}^{\al(i)}\circ p_i^{\rho(\al(i))}$.  Clearly
$pro(f_\bul)=(\bar h_{\rho(j)}^j)_{j\in J}$.
\qed

\begin{corollary}
\label{approxiffpro}
Suppose $I,J\cong\Z_+$ and $f_\bul:A_\bul\ra B_\bul$ is a morphism.
  Then $f_\bul$ is an approximate
isomorphism iff $pro(f_\bul)$ is a pro-isomorphism.
\end{corollary}
\proof  By Lemma \ref{isoiffapproxinvertible}, $f_\bul$ is an approximate isomorphism
iff it represents an invertible element of $Approx$, and by 
Lemma \ref{faithful} this is equivalent to saying that $pro(f_\bul)$ is
invertible in {\it Pro-Abelian}.
\qed

\subsection{Recognizing  groups of type $FP_n$}
\label{recognitionsection}

The main result in this section is Theorem \ref{finitenesstheorem},
which gives a characterization of groups $G$ of
type $FP_n$ in terms of  nested families
of $G$-chain complexes, and Lemma \ref{finlem} which 
relates the cohomology of $G$
with the corresponding cohomology of the $G$-chain complexes.
A related characterization of groups of type $FP_n$
appears in \cite{brown1}.  We will apply
Theorem \ref{finitenesstheorem} and Lemma \ref{finlem} 
in section \ref{mainargument}.    

Suppose for $i=0,\ldots,N$ we have an augmented 
chain complex $A_*(i)$ of projective
$\Z G$-modules, and for $i=1,\ldots,N$ we have an augmentation
preserving $G$-equivariant chain map $a_i:A_*(i-1)\ra A_*(i)$ which induces 
zero on reduced homology in dimensions $< n$. 
Let $G$ be a group of type $FP_k$, and let
$$0\larrow \Z\larrow P_0\larrow\ldots\larrow P_k$$
be a partial resolution $P_*$ of $\Z$ by finitely generated
projective $\Z G$-modules. 
We assume that $k\le n \le N$. 

\begin{lemma}
\label{finlem}
Under  the above conditions we have:

1. There is an augmentation 
preserving $G$-equivariant chain mapping $P_*\ra A_*(n)$.

2. If $k< n$ and $j_i: P_*\ra A_*(0)$
are augmentation preserving $G$-equivariant chain mappings
for $i=1,\,2$, then the compositions $P_*\stackrel{j_i}{\ra}
A_*(0)\ra A_*(k)$ are $G$-equivariantly chain homotopic.

3. \sloppy{Suppose $k< n$ and $f:P_*\ra A_*(0)$ is an augmentation 
preserving $G$-equivariant chain mapping.  Then 
for any $\Z G$-module $M$, the map 
$$
H^i(f):H^i(A_*(0);M)\ra H^i(P_*;M)$$
carries the image $Im(H^i(A_*(n);M)\ra H^i(A_*(0);M))$ isomorphically
onto $H^i(P_*;M)$ for $i=0,\ldots k-1$.  The map 
$$
H_i(f):H_i(P_*;M)\ra H_i(A_*(n);M)$$ 
is an  isomorphism onto 
the image of  $H_i(A_*(0);M)\ra H_i(A_*(n);M)$ for $i=0,\ldots k-1$. 
The map 
$$H_k(f):H_k(P_*;M)\ra H_k(A_*(n);M)$$
is onto the image of $H_k(A_*(0);M)\ra H_k(A_*(n);M)$.} 
\end{lemma}
{\em Proof of 1.} We start with the diagram 
$$
\begin{array}{cc}
P_0 & \ \\
\downarrow & \ \\
\Z & \leftarrow A_0(0) . 
\end{array}
$$
Then projectivity of $P_0$ implies that we can complete this to a commutative 
diagram by a $\Z G$-morphism $f_0: P_0\to A_0(0)$.  
Assume inductively that we have constructed
a $G$-equivariant augmentation preserving chain mapping 
$f_j:[P_*]_j\ra A_*(i)$.
Then the image of the composition 
$P_{i+1}\stackrel{\D}{\ra}P_j\stackrel{f_j}{\ra}A_j(j)\ra A_j(j+1)$ 
is contained in the image of 
$A_{j+1}(j+1)\stackrel{\D}{\ra}A_j(j+1)$
since $a_{j+1}$ induces zero on reduced homology. So projectivity of 
$P_{j+1}$ allows us to extend $f_j$ to a $G$-equivariant chain 
mapping $f_{j+1}:[P_*]_{j+1}\ra A_*(j+1)$.

{\em Proof of 2.} Similar to the proof of 1: use induction and
projectivity of the $P_\ell$'s.

{\em Proof of 3.} Let 
$\rho_*:[A_*(n)]_k\ra P_*$ be a $G$-equivariant
chain mapping constructed using the fact that $H_i(P_*)=\{0\}$ for $i<k$.
Consider the compositions
$$
\al_{k-1}: [P_*]_{k-1}\stackrel{f_*}{\ra}[A_*(0)]_{k-1}\ra 
[A_*(n)]_{k-1}\stackrel{\rho_*}{\ra}P_*
$$
and
$$
\be_k: [A_*(0)]_k\ra[A_*(n)]_k\stackrel{\rho_*}{\ra}
[P_*]_k\stackrel{f_*}{\ra}[A_*(0)]_k \ra A_*(n).
$$
Both are ($G$-equivariantly) chain homotopic to the inclusions; 
the first one since $P_*$ is a partial resolution, and the second by applying  
assertion 2 to the chain mapping $[A_*(0)]_k\ra A_*(0)$.  
Assertion 3 follows immediately from this.
\qed

We note that this lemma did not require any finiteness assumptions on the 
$\Z G$-modules $A_i(j)$. Suppose now that the group $G$ satisfies the above 
assumptions and  let $G\acts X$ be a free simplicial action
on a uniformly $(n-1)$-acyclic locally finite metric simplicial complex $X$, 
$k\le n-1$. 
Then by part 1 of the previous lemma we have a $G$-equivariant augmentation-preserving 
chain mapping $f:P_*\ra C_*(X)$. Let $K\subset X$ be the support of the 
image of $f$. It is clear that $K$ is  $G$-invariant and $K/G$ is compact. 
As a   corollary of the proof of the previous lemma, we get: 

\begin{corollary}
\label{co} 
Under the above assumptions the direct system of reduced homology groups 
$\{\t H_i(N_R(K))\}_{R\ge 0}$ is approximately zero for each $i<k$.  
\end{corollary}
\proof Given $R>0$ we consider the system of chain complexes 
$A_*(0):= C_*(N_R(K))$, $A_*(1)=A_*(2)=... =A_*(N)= C_*(X)$. 
The mapping $[A_*(0)]_k \stackrel{\be_k}{\to} A_*(N)= C_*(X)$ from the proof 
of Lemma \ref{finlem} is chain homotopic  to 
the inclusion via a $G$-equivariant homotopy $h_R$. On the other hand, this 
map factors through $P_*$, hence it induces zero mapping of the reduced 
homology groups
$$
\t H_i(N_R(K))\stackrel{0}{\to} \t H_i(Support(Im(\be_k))), \ i<k. 
$$ 
The support of $Im(h_R)$ is contained in  $N_{R'}(K)$ for some $R'<\infty$, 
since  $h_R$ is $G$-equivariant. Hence the inclusion $N_R(K)\to  N_{R'}(K)$ 
induces zero map of $\t H_i(\cdot )$ for $i<k$. \qed 

Before stating the next corollary, we recall the following
fact:

\begin{lemma}
\label{cptsuppiso}(See \cite{Brown}.)
Let $G\acts X$ be  a discrete, free, cocompact action of a group on a simplicial
complex.  Then the complex of compactly supported simplicial cochains
$C^*_c(X)$ is canonically isomorphic to the complex $Hom_{\Z G}(C_*(X);\Z G)$; 
in particular, the compactly supported cohomology of $X$ is canonically isomorphic to $H^*(X/G;\Z G)$.
\end{lemma}

In the  next corollary we assume that $G$, $P_*$, $X$, $f$, $K$ are as above, 
in particular, $X$ is a  uniformly $(n-1)$-acyclic locally finite metric simplicial 
complex, $k\le n-1$, $P_k \to ... \to P_0 \to \Z\to 0$ is a resolution 
by finitely generated projective $\Z G$ modules.   

\begin{corollary}
\label{retractionres}
1. For any local coefficient 
system ($\Z G$-module)  $M$ the  family of maps 
$$H^i(N_R(K)/G;M)\stackrel{f_R^i}{\ra} H^i(P_*;M)$$
defines a morphism  between
the inverse system $\{H^i(N_R(K)/G;M)\}_{R\geq 0}$
and the constant inverse system $\{ H^i(P_*;M)\}_{R\geq 0}$
which  is an approximate isomorphism 
when $0\leq i<k$. 

2.  The map
$$
H_c^i(N_R(K))\simeq H^i(N_R(K)/G;\Z G)\stackrel{f_R^i}{\lra} H^i(P_*;\Z G)$$
is an approximate isomorphism 
when $0\leq i<k$. 

3.  The $\Z G$-chain map 
$$
f_{R,*}:P_*\ra C_*(N_R(K))
$$
induces a homomorphism of homology groups
$$
f_{R,i}:\tilde H_i(P_*,\Z G)\ra \tilde H_i(N_R(K))
$$
which is an approximate isomorphism for $0\leq i<k$. 
\end{corollary}
\proof 1. According to Corollary \ref{co}  the direct system of 
reduced homology groups $\{\t{H}_i(N_R(K))\}$ is pro-zero for each 
$i<k$. Thus for $N>k$ we have a sequence of integers 
$R_0=0 < R_1 < R_2 <  ...< R_N$ so that the maps 
$$
\t{H}_i(N_{R_j}(K))\to \t{H}_i(N_{R_{j+1}}(K))
$$ 
are zero for each $j<N, i<k$. We now apply Lemma \ref{finlem} where 
$A_*(j):= C_*(N_{R_j}(K))$.                   

2. This follows from part 1 and Lemma \ref{cptsuppiso}.

3. Note that $\tilde H_i(P_*,\Z G)\simeq\{0\}$ for $i<k$; this 
follows directly from the definition of a group of type $FP_k$.
Thus the assertion follows from Corollary \ref{co}. \qed

There is also an analog of Corollary \ref{retractionres}
which does not require a group action:

\begin{lemma}
\label{retractionmap}
Let $X$ and $Y$ be bounded geometry metric simplicial complexes,
where $Y$ is uniformly $(k-1)$-acyclic and 
$X$ is uniformly $k$-acyclic.  Suppose $C_*(Y)\stackrel{f}{\ra} C_*(X)$
is a uniformly proper chain mapping, and $K\defeq Support(Im(f))\subset X$.
Then 

1. The induced map on cohomology
$$H^i_c(f):H^i_c(N_R(K))\ra H^i_c(Y)$$
defines a morphism  between
the inverse system $\{H^i_c(N_R(K))\}_{R\geq 0}$
and the constant inverse system $\{ H^i_c(Y)\}_{R\geq 0}$ 
which is an approximate isomorphism for $0\leq i<k$,
and an approximate monomorphism for $i=k$.

2. The approximate isomorphism approximately
respects support in the following sense.  There is a function $\zeta:\N\ra\N$
so that if  $i<k$, $S\subset Y$ is a subcomplex, 
$T\defeq Support(f_*(C_*(S)))\subset X$
is the corresponding subcomplex of $X$, and
$\al\in Im(H^i_c(Y,\ol{Y-S})\ra H^i_c(Y))$, 
then $\al$ belongs to the image of the composition
$$
H^i_c(N_R(K),\ol{N_R(K)-N_{\zeta(R)}(T)})
\ra H^i_c(N_R(K))\stackrel{H^i_c(f)}{\RA}H^i_c(Y).
$$

3. The induced map 
$$
\tilde H_i(f):\{0\}\simeq \t H_i(Y)\ra \t H_i(N_R(K))
$$
is an approximate isomorphism for $0\leq i<k$.
\end{lemma}
\proof
Since $f$ is uniformly proper, using  
the uniform $(k-1)$-acyclicity of $Y$ and
uniform $k$-acyclicity of $X$, we can 
construct a direct system $\{\rho_R\}$ of uniformly proper chain mappings
 between the truncated chain
complexes
$$[0\leftarrow C_0(N_R(K))\leftarrow\ldots\leftarrow
C_k(N_R(K))]\stackrel{\rho_R}{\ra}
[0\leftarrow C_0(Y)\leftarrow\ldots\leftarrow
C_k(Y)]$$
so that the compositions $f\circ\rho_R$ are 
chain homotopic to the inclusions 
$$[0\leftarrow C_0(N_R(K))\leftarrow\ldots\leftarrow
C_k(N_R(K))]$$
$$\ra [0\leftarrow C_0(N_{R'}(K))\leftarrow\ldots\leftarrow
C_k(N_{R'}(K))]$$
(for $R'=\om(R)$) via chain homotopies of bounded support. 
Moreover the restriction of the composition $\rho_R\circ f$ 
to the $(k-1)$-truncated chain complexes is chain
homotopic to the identity via a chain homotopy with bounded support.

We first prove that the morphism of inverse systems defined by 
$$
H^i_c(f):H^i_c(N_R(K))\ra H^i_c(Y)
$$
is an approximate monomorphism.  Suppose 
$$\al\in
Ker(H^i_c(f):H^i_c(N_{R'}(K))\ra H^i_c(Y))$$ 
where $R'=\om(R)$. 
Then $H^i(f\circ\rho_{R'})(\al)=0$.  But
the restriction of $H^i(f\circ\rho_{R'})(\al)$
to $N_R(K)$ is cohomologous to the 
restriction of $\al$ to $N_R(K)$.

Since the restriction of the composition $\rho_R\circ f$ 
to the $(k-1)$-truncated chain complex $[C_*(Y)]_{k-1}$ 
is chain homotopic to the identity, it follows that 
$$H^i_c(f):H^i_c(N_R(K))\ra H^i_c(Y)$$
is an epimorphism for $R\geq 0$ and $i<k$.

Part 2 of the lemma follows immediately
from the uniform properness of $\rho_R$ and the coarse
Lipschitz property of the chain homotopies constructed above.

We omit the proof of part 3 as it is similar to that of part 2.
\qed

Note that in the above discussion 
we used finiteness  assumptions on the group $G$ 
to make conclusions about (co)homology of families of $G$-invariant 
chain complexes.  
Our next goal is to use existence of a family of chain complexes $A_*(i)$ 
of {\em finitely generated} 
projective $\Z G$ modules as in Lemma \ref{finlem} to establish finiteness 
properties of the group $G$ (Theorem \ref{finitenesstheorem}). 
We begin with a homotopy-theoretic analog of Theorem 
\ref{finitenesstheorem}. 

\begin{proposition}
Let $G$ be a group, and let $X(0)\stackrel{a_1}{\ra} X(1)\stackrel{a_2}{\ra}
 \ldots\stackrel{a_{n+1}}{\ra} X(n+1)$
be a diagram of free, simplicial $G$-complexes
where $X(i)/G$ is  compact for $i=0,\ldots n+1$.  If 
the maps $a_i$ are $n$-connected for each $i$, then
there is an $(n+1)$-dimensional free, simplicial $G$-complex
$Y$ where $Y/G$ is compact and $Y$ is $n$-connected.
\end{proposition}

\proof  We build $Y$ inductively as follows.  Start
with $Y_0=G$ where $G$ acts  on $Y_0$ by left translation,
 and let $j_0:Y_0\ra X(0)$ be any
$G$-equivariant simplicial map.  Inductively apply Lemma \ref{homotopyinduction}
below to the composition  $Y_i\stackrel{j_i}{\ra} X(i)\ra X(i+1)$
to obtain $Y_{i+1}$ and a simplicial $G$-map $j_{i+1}:Y_{i+1}\ra X(i+1)$.
Set $Y\defeq Y_{n+1}$.
\qed

\begin{lemma}
\label{homotopyinduction}
Let $Z$ and $A$ be locally finite simplicial
complexes with free cocompact simplicial
$G$-actions, where $dim(Z)=k$, and $Z$ is $(k-1)$-connected.
 Let $j:Z\ra A$, be a null-homotopic
$G$-equivariant  simplicial map.
Then we may construct a $k$-connected simplicial
$G$-complex $Z'$
by attaching (equivariantly) finitely many $G$-orbits 
of   simplicial \footnote{A simplicial cell
is a simplicial complex PL-homeomorphic to
a single simplex.} $(k+1)$-cells to $Z$, and a $G$-map $j':Z'\ra A$
extending $j$.
\end{lemma}
\proof
By replacing $A$ with the mapping cylinder of $j$,
we may assume that $Z$ is a subcomplex of $A$ and
$j$ is the inclusion map.  Let $A_k$ denote the
$k$-skeleton of $A$.   Since $Z$ is $(k-1)$-connected,
after subdividing $A_k$ if necessary, 
we may construct a $G$-equivariant  simplicial retraction
$r:A_k\ra Z$.  
For every $(k+1)$-simplex $c$ in $A$,
we attach a  simplicial $(k+1)$-cell $c'$ to $Z$ using the
composition of the attaching map of $c$ with 
the retraction $r$.  It is clear that we may do this
$G$-equivariantly, and there will be only finitely
many $G$-orbits of $(k+1)$-cells attached.  We denote
the resulting simplicial complex by $Z'$, and note that
the inclusion $j:Z\ra A$ clearly
extends (after subdivision of $Z'$) to an equivariant 
simplicial map $j':Z'\ra A$.
 
 We now claim that $Z'$ is $k$-connected.  Since
 we built $Z'$ from $Z$ by attaching $(k+1)$-cells,
 it suffices to show that $\pi_k(Z)\ra\pi_k(Z')$ is trivial.
  If $\si:S^k\ra Z$ is a simplicial map for
  some triangulation of $S^k$, we get a simplicial null-homotopy 
$\tau:D^{k+1}\ra A$ extending $\si$. 
Let $D^{k+1}_k$ denote the $k$-skeleton of $D^{k+1}$.
The composition $D^{k+1}_k\stackrel{\tau}{\ra}A\stackrel{r}{\ra}
 Z\ra Z'$
extends over each $(k+1)$-simplex $\De$ of $D^{k+1}$, since
$\tau\restr_{\De}:\De\ra A$ is either an 
embedding, in which case $r\circ\tau\restr_{\D\De}:\D\De\ra Z'$
is null homotopic by the construction of $Z'$, or 
$\tau\restr_{\De}:\De\ra A$ has image contained in a $k$-simplex
of $A$, and the composition 
$$\D\De\stackrel{\tau}{\ra}A\stackrel{r}{\ra}Z$$
is already null-homotopic.  Hence the composition
$S^k \stackrel{\si}{\ra}Z\hookrightarrow Z'$ is null-homotopic.
\qed

The next lemma is a homological analog of Lemma \ref{homotopyinduction}
which provides the inductive step in the proof of 
Theorem \ref{finitenesstheorem}.

\begin{lemma}
Let $G$ be a group.   Suppose 
$0\larrow \Z\stackrel{\eps}{\larrow}P_0\larrow\ldots\larrow P_k$
is a partial resolution by finitely generated projective
$\Z G$-modules, and $\Z\stackrel{\eps}{\larrow}A_0\larrow\ldots\larrow A_{k+1}$
is an augmented chain complex of finitely generated projective
$\Z G$-modules.  Let $j:P_*\ra A_*$ be an augmentation preserving 
chain mapping which induces zero on homology groups\footnote{We
declare that $H_k(P_*)\defeq Z_k(P_*)$.}.  Then 
we may extend $P_*$ to a partial resolution $P_*'$: 
$$
0 \larrow \Z\stackrel{\eps}{\larrow}P_0\larrow\ldots\larrow P_k\larrow P_{k+1}
$$
where $P_{k+1}$ is finitely generated free, and $j$ extends
to a chain mapping $j':P_*'\ra A_*$.
\end{lemma}
\proof
By replacing $A_*$ with the algebraic mapping
cylinder of $j$, we may assume that $P_*$ is embedded
as a subcomplex of $A_*$, $j$ is  the inclusion,
and for $i=0,\ldots,k$, the chain group 
$A_k$ splits as a direct sum of $\Z G$-modules
$A_i=P_i\oplus Q_i$ where $Q_i$ is finitely generated
and projective.  Applying the projectivity
of $Q_i$, we construct a chain retraction 
from the $k$-truncation
$[A_*]_k$ of $A_*$ to $P_*$.  Choose a finite set
of generators $a_1,\ldots,a_\ell$ for the $\Z G$-module
$A_{k+1}$.  We let $P_{k+1}$ be the free module
of rank $\ell$, with basis $a_1',\ldots,a_\ell'$, and
define the boundary operator $\D:P_{k+1}\ra P_k$
by the formula $\D(a_i')=r(\D(a_i))$.  To see
that $H_k(P_*')=0$, pick a $k$-cycle $\si\in Z_k(P_*)$.
We have $\si=\D\tau$ for some $\tau=\sum c_ia_i\in A_{k+1}$.
Then $\si=r(\D\tau)=\sum c_ir(\D a_i)=\sum c_i\D a_i'$;
so $\si$ is null-homologous in $P_*'$. The extension mapping 
$j':P_*'\ra A_*$ is defined  by $a_i'\mapsto a_i, 1\le i\le \ell$. 
\qed

\begin{theorem}
\label{finitenesstheorem}
Suppose for $i=0,\ldots,N$ we have an augmented 
chain complex $A_*(i)$ of finitely generated projective
$\Z G$-modules, and for $i=1,\ldots,N$ we have an augmentation
preserving $G$-equivariant chain map $a_i:A_*(i-1)\ra A_*(i)$ 
which induces zero on reduced homology in dimensions $\leq n \le N$.   

Then there is a partial resolution 
$$0\larrow \Z\larrow F_0\larrow\ldots\larrow F_n$$
of finitely generated free $\Z G$-modules, and 
a $G$-equivariant chain mapping $f:F_*\ra A(n)$.  In particular,
$G$ is a group of type $FP_n$.
\end{theorem}
\proof Define $F_0$ to be the group ring $\Z G$, with the 
usual augmentation $\Z\larrow \Z G$.  Then
construct $F_i$ and a chain map $F_i \to A_i(i)$ by 
applying the previous lemma inductively.
\qed

\begin{corollary}
Suppose that $G\acts X$ is a free simplicial action of a group $G$ 
on a metric simplicial complex $X$. Suppose that we have a system of 
(nonempty) $G$-invariant simplicial subcomplexes $X(0)\subset X(1) \subset ...
\subset X(N)$ so that:

(a) $X(i)/G$ is compact for each $i$,

(b) The induced mappings $\t {H}_i(X(k))\to \t {H}_i(X(k+1))$ 
are zero for each $i\le n\le N$ and $0\le k < N$. 

Then the group $G$ is of type $FP_n$.  
\end{corollary}
\proof Apply Theorem \ref{finitenesstheorem}  to $A_*(i):= C_*(X(i))$. \qed  

Note that the above corollary is the converse to Corollary \ref{co}. Thus 

\begin{corollary}
Suppose that $G\acts X$ is a group action on a uniformly acyclic bounded 
geometry metric simplicial complex, $K:= G(\star)$, where $\star\in X$. 
Then $G$ is of type $FP$ if and only 
if the the direct system of reduced homology groups $\{\t{H}_*(N_R(K))\}$ 
is pro-zero.   
\end{corollary}

Combining  Theorem \ref{finitenesstheorem} and Lemma \ref{finlem} we get: 

\begin{corollary}
\label{finitenesscor}
Suppose for $i=0,\ldots,2n+1$ we have an augmented 
chain complex $A_*(i)$ of finitely generated projective
$\Z G$-modules, and for $i=1,\ldots,2n+1$ we have  augmentation
preserving $G$-equivariant chain maps $a_i:A_*(i-1)\ra A_*(i)$ which induce zero
on reduced homology in dimensions $\leq n$.  Then: 

1.  There is a partial resolution $F_*$:
$$0\larrow\Z\larrow F_0\larrow\ldots\larrow F_n$$
by finitely generated free $\Z G$-modules and a $G$-equivariant
chain mapping $f_*:F_*\ra A_*(n)$.  In particular
$G$ is of type $FP_n$.

2.  For any $\Z G$-module $M$, the map 
$H^i(f):H^i(A_*(n);M)\ra H^i(F_*;M)$ carries
the image $Im(H^i(A(2n);M)\ra H^i(A(n);M))$ isomorphically
onto $H^i(F_*;M)$ for $i=0,\ldots n-1$.

3. The map 
$H_i(f):H_i(P_*;M)\ra H_i(A_*(2n);M)$ is an isomorphism onto 
the image of  $H_i(A_*(n);M)\ra H_i(A_*(2n);M)$. 
\end{corollary}

We now discuss a relative version of Corollaries \ref{retractionres} and 
\ref{finitenesscor}.  
Let $X$ be a uniformly acyclic bounded geometry metric simplicial complex, 
and $G$ be a group  
acting freely simplicially on $X$; thus $G$ has 
finite cohomological dimension since $X$ is acyclic and $dim(X)<\infty$. 
Let $K\subset X$ be a 
$G$-invariant subcomplex  so that $K/G$ is compact; and let 
$\{C_\al\}_{\al\in I}$ be the deep  components of $X-K$. 
Define  $Y_R:= \ol{X- N_R(K)}$, $Y_{\al,R}:= C_\al \cap Y_R$. 
We will assume that the system
$$
\{\t H_j(Y_{\al,R})\}_{R\geq 0}
$$
is approximately zero for each $j,\al$. In particular, 
$\{\t H_0(Y_{\al,R})\}_{R\geq 0}$ is approximately zero, which implies that
each $C_\al$ 
is stable. Let $H_\al$ denote the stabilizer of $C_\al$ in $G$. 
Choose a set of representatives  $C_{\al_1},\ldots,C_{\al_k}$
from the  $G$-orbits in the collection $\{C_\al\}$.   For notational
simplicity we relabel $\al_1,\ldots,\al_k$ as $1,\ldots,k$.  
Let $H_i=H_{\al_i}$ be  the stabilizer of $C_i=C_{\al_i}$.    
This defines a group pair $(G; H_1,..., H_k)$. 
Let $P_*$ be a finite length projective resolution of $\Z$ by 
$\Z G$-modules, and for each $i=1,\ldots,k$, we choose 
 a finite length projective resolution of $\Z$ by 
$\Z H_i$-modules $Q_*(i)$.    Using the construction described in 
section \ref{groupprelim} (see the discussion of the group pairs) 
we convert this data to a pair
$(C_*,D_*)$ of finite length projective resolutions (consisting
of $\Z G$-modules).    We recall that $D_*$ decomposes in a natural
way as a direct sum $\oplus_\al D_*(\al)$ where each $D(\al)$ is a
resolution of $\Z$ by projective $\Z H_\al$-modules.   Now
construct a $\Z H_i$-chain mapping $C_*(Y_{\al_i,0})\ra D_*(\al_i)$
using the acyclicity of $D_*(\al_i)$.   We then extend this $G$-equivariantly
to a mapping $C_*(Y_0)\ra D_*$, and then to a $\Z G$-chain mapping 
$\rho_0:(C_*(X), C_*(Y_0))\to (C_*,D_*)$.    By restriction, this defines 
a morphism of inverse systems $\rho_R:(C_*(X), C_*(Y_R))\to (C_*,D_*)$. 

\begin{lemma}
\label{relretractionres}
The mapping $\rho_\bul$ induces approximate isomorphisms between
relative (co)homology with local coefficients:
$$
H^*(G,\{H_i\}; M)\ra   H^*(C_*(X),C_*(Y_R);M)\simeq H^*(X/G,Y_R/G;M) $$ 
$$H_*(X/G,Y_R/G;M)\simeq H_*(C_*(X),C_*(Y_R);M) \ra H_*(G,\{H_i\};M)$$
for any $\Z G$-module $M$.
\end{lemma}
\proof We will prove the lemma by showing that the maps $\rho_R$ form
an ``approximate chain homotopy equivalence'' in an appropriate sense.

For each $i$ we construct  a $\Z H_i$-chain mapping $D_*(i)\ra C_*(Y_{i,R})$
using part 1 of Lemma \ref{finlem} and  the fact that 
$$
\{\t H_j(Y_{\al,R})\}_{R\geq 0}
$$
is an approximately zero system. We  then extend these to 
$\Z G$-chain mappings
$$
f_R: (C_*, D_*) \to (C_*(X), C_*(Y_R)).
$$
Using part 2 of Lemma \ref{finlem}, we can
actually choose the  mappings $f_R$ so that they form a compatible system 
chain mappings up to chain-homotopy. 
The composition
$$
\rho_R \circ f_R : (C_*,D_*)\to (C_*,D_*)
$$
is  $\Z G$-chain mapping, hence it is chain-homotopic to the
identity. The composition
$$
f_R \circ \rho_R: C_*(X,Y_{R})\to C_*(X,Y_{R})
$$
need not be chain homotopic to the identity, but it becomes chain
homotopic to the projection map when
precomposed with the restriction $C_*(X,Y_{R'})\ra C_*(X,Y_R)$
where $R'\geq R$ is suitably chosen (by again using part 2 of
Lemma \ref{finlem} and the fact that  
$$
\{\t H_j(Y_{\al,R})\}_{R\geq 0}
$$
is an approximately zero system).     This clearly implies the induced
homorphisms on (co)homology are approximate isomorphisms.\qed

\subsection{Coarse Poincare duality}
\label{cpd}

We now introduce a class of metric simplicial complexes
which satisfy coarse versions of Poincare and Alexander duality, see
Theorems \ref{Coarse Poincare duality},  \ref{aduality},  
\ref{adualitymapping}.

From now on we will adopt the  convention of extending each (co)chain 
complex indexed by the nonnegative integers to a complex
indexed by the integers by setting the remaining groups equal to zero.
 So  
for each (co)chain complex $\{C_i, i\ge 0\}$ we get the (co)homology 
groups $H_i(C_*), H^i(C_*)$ 
defined for $i< 0$.

\begin{definition}[Coarse Poincar\'e duality spaces]
\label{coarsepdndef}
\sloppy{A {\em coarse Poincar\'e duality space of formal
dimension} $n$ is 
 a bounded geometry metric simplicial complex $X$ so that
$C_*(X)$ is uniformly acyclic, and there is a constant $D_0$ 
and  chain mappings 
$$
C_*(X)\stackrel{\bar P}{\ra}C_c^{n-*}(X)\stackrel{P}{\ra}C_*(X)
$$
so that}

1. $P$ and $\bar P$ have displacement $\leq D_0$ (see section \ref{geomprelim}
for the definition of displacement).

2.   $\bar P\circ P$ and $P\circ\bar P$ are
chain homotopic to the identity by $D_0$-Lipschitz\footnote{See 
section \ref{geomprelim}.} chain homotopies 
$\Phi:C_*(X)\ra C_{*+1}(X)$,
$\bar\Phi:C^*_c(X)\ra C^{*-1}_c(X)$.

We will often refer to coarse Poincare duality spaces of formal
dimension $n$ as {\em coarse $PD(n)$ spaces.}  
Throughout the paper we will reserve the letter $D_0$ for the constant which appears
in the definition of a coarse $PD(n)$ space; we let $D\defeq D_0+1$.  
\end{definition}

Note that for each coarse $PD(n)$ space $X$ we have
$$
H^*_c(X)\simeq H_{n-*}(X)\simeq H_{n-*}(\R^n)\simeq H^*_c(\R^n).$$ 
We will not need
the bounded geometry and uniform acyclicity conditions until
Theorem \ref{adualitymapping}.  Later in the paper
we will consider simplicial actions on coarse $PD(n)$ spaces,
and we will assume implicitly that the actions commute with the
operators $\bar P$ and $P$, up to  chain homotopy with uniformly bounded
Lipschitz constants.

The next lemma gives  important examples of coarse $PD(n)$
spaces: 
  
\begin{lemma}
\label{thesearepdn}
The following are coarse $PD(n)$ spaces:

1. An acyclic metric simplicial complex $X$ which admits
a free, simplicial, cocompact action by a $PD(n)$ group. 

2.  An $n$-dimensional, bounded geometry metric simplicial
complex $X$, with an augmentation $\al:C^n_c(X)\ra\Z$
for the compactly supported simplicial cochain complex,
so that $(C^*_c(X),\al)$ is uniformly acyclic (see section \ref{geomprelim} for definitions). 

3. A uniformly acyclic, bounded geometry metric simplicial complex $X$
which is a topological $n$-manifold.
\end{lemma}
{\em Proof of 1.}  Let $\pres{n}$ be a resolution of $\Z$
by finitely generated projective $\Z G$-modules.  $X$ is 
acyclic, so we have $\Z G$-chain homotopy equivalences 
 $P_*\stackrel{\al}{\simeq} C_*(X)$
and $Hom(P_*,\Z G)\simeq C^*_c(X)$ where $\al$ is augmentation
preserving.  
Hence to construct the two chain equivalences needed in 
Definition \ref{coarsepdndef}, it suffices to 
construct a $\Z G$-chain homotopy equivalence $p:P_*\ra Hom(P_{n-*},\Z G)$
of $\Z G$-modules (since the operators are $G$-equivariant
conditions 1 and 2 of Definition \ref{coarsepdndef} will
be satisfied automatically).  For this, see \cite[p. 221]{Brown}.

{\em Proof of 2.}  We construct a chain mapping $P:C_*(X)\ra C^{n-*}_c(X)$
as follows.   We first map each vertex $v$ of $X$ to an $n$-cocycle
$\be\in C^n_c(X,\ol{X-B(v,R_0)})$ which maps to $1$ under the 
augmentation $\al$,  (such a $\be$ exists by the
uniform acyclicity of $(C^*_c(X),\al)$), and extend this to
a homomorphism $C_0(X)\ra C^n_c(X)$.  By the uniform acyclicity
of $(C^*_c(X),\al)$ we can extend this to a chain mapping $P$.
By similar reasoning we obtain a chain homotopy inverse $\bar P$, and
construct chain homotopies $\bar P\circ P\sim id$ and
$P\circ \bar P\sim id$.

{\em Proof of 3.}  $X$ is acyclic, and therefore orientable.  
An orientation of $X$ determines an augmentation $\al:C^n_c(X)\ra\Z$.
The uniform acyclicity of $X$ together with ordinary Poincare duality implies
that $(C^*_c(X),\al)$ is uniformly acyclic.  So 3 follows from 2.

We remark that if $G\acts X$ is a free simplicial action then
these constructions can be made $G$-invariant. \qed

\medskip
When $K\subset X$ is a (nonempty) subcomplex we will consider the direct system
of tubular neighborhoods $\{N_R(K)\}_{R\ge 0}$ of $K$ and the inverse system 
of the closures of their complements 
$$
\{Y_R:= \ol{X- N_R(K)}\}_{R\ge 0}.
$$ 
We get four  inverse and four direct systems of (co)homology groups:
$$
\{H^k_c(N_{R}(K))\},  \,\{H_{j}(X,Y_R)\},\,  \{H^k_c(X, N_{R}(K))\}, \, \{H_{j}(Y_{R})\}
$$
$$
\{H^k_c(Y_R)\} ,\,  \{H_{j}(X, N_{R}(K))\},  \,\{H^k_c(X, Y_{R})\},  \,\{H_{j}(N_{R}(K))\}
$$
 with the usual restriction and projection homomorphisms.
Note that by excision, we have  isomorphisms
$$
H_{j}(X,Y_R)\simeq H_{j}(N_R(K),\D N_R(K)), \hbox{~~etc}. 
$$

Extension by zero defines a group homomorphism
$C^k_c(N_{R+D}(K))\stackrel{ext}{\subset}C^k_c(X)$. When we compose this with 
$$
C^k_c(X)\stackrel{P}{\ra} 
C_{n-k}(X)\stackrel{proj}{\ra} C_{n-k}(X,Y_R) 
$$
we get a well-defined induced homomorphism
$$
P_{R+D}:H^k_c(N_{R+D}(K))\ra H_{n-k}(X,Y_R)
$$
where $D$ is as in Definition \ref{coarsepdndef}.
We get, in a similar fashion, homomorphisms
\begin{equation}
\label{e1}
H^k_c(N_{R+D}(K)) \stackrel{P_{R+D}}{\lra} H_{n-k}(X,Y_R) \stackrel{\bar P_{R}}{\lra}
H^k_c(N_{R-D}(K)) 
\end{equation}
\begin{equation}
\label{e2}
H^k_c(Y_R) \stackrel{P_{R}}{\lra} H_{n-k}(X, N_{R+D}(K)) \stackrel{\bar P_{R+D}}{\lra}
H^k_c(Y_{R+2D})
\end{equation}
\begin{equation}
\label{e3}
H^k_c(X, N_{R+D}(K)) \stackrel{P_{R+D}}{\lra} H_{n-k}(Y_{R}) \stackrel{\bar P_{R}}{\lra}
H^k_c(X, N_{R-D}(K)) 
\end{equation}
\begin{equation}
\label{e4}
H^k_c(X, Y_{R}) \stackrel{P_{R}}{\lra} H_{n-k}(N_{R+D}(K)) \stackrel{\bar P_{R+D}}{\lra}
H^k_c(X, Y_{R+2D}) 
\end{equation}
Note that the homomorphisms in (\ref{e1}), (\ref{e3})  determine $\al$-morphisms between inverse systems 
and the homomorphisms in (\ref{e2}), (\ref{e4})  determine $\be$-morphisms between direct systems, 
where $\al(R)=R-D$, $\be(R)=R+D$ (see section \ref{algprelim} for definitions).
 These operators inherit the bounded displacement property
of $P$ and $\bar P$, see condition 1 of Definition \ref{coarsepdndef}.
We let $\om(R):= R+2D$, where $D$ is the constant from Definition \ref{coarsepdndef}.

\begin{theorem}
[Coarse Poincare duality]
\label{Coarse Poincare duality}
Let $X$ be a coarse $PD(n)$ space, $K\subset X$ be a subcomplex as above. Then the morphisms 
$P_\bullet, \bar P_\bullet$ in (\ref{e1}), (\ref{e3}) are $(\om,\om)$-approximate isomorphisms of inverse systems 
and the morphisms 
$P_\bullet, \bar P_\bullet$ in (\ref{e2}), (\ref{e4}) are $(\om,\om)$-approximate isomorphisms of direct systems
(see section \ref{algprelim}). 
In particular, if $X\ne N_{R_0}(K)$ for any $R_0$ then the inverse systems $\{H^n_c(N_R(K))\}_{R\ge 0}$ and 
$\{H_n(Y_R)\}_{R\ge 0}$  are approximately zero. 
\end{theorem}
\proof We will verify the assertion for the homomorphism $P_\bullet$ in (\ref{e1}) and leave the rest to the reader. 
We first check that $P_\bullet$ is an $\om$-approximate 
monomorphism.  Let 
$$
\xi\in Z^*_c(N_{R+2D}(K))$$
be a cocycle representing an element  $[\xi]\in Ker(P_{R+2D})$, and let
$\xi_1\in C^*_c(X)$ be the extension of $\xi$ by zero.
Then we have 
$$P(\xi_1)=\D\eta +\zeta$$
where $\eta\in C_{n-*}(X)$ and $\zeta\in C_{n-*}(\ol{X-N_{R+D}(K)})$.
Applying $\bar P$ and the chain homotopy $\Phi$,
we get 
$$
\de\Phi(\xi_1)+\Phi\de(\xi_1)=\bar P\circ P(\xi_1)-\xi_1
=\bar P(\D\eta+\zeta)-\xi_1
$$
so
$$
\xi_1=\de\bar P(\eta)+\bar P(\zeta)-\de\Phi(\xi_1)-\Phi\de(\xi_1).
$$
The second and fourth terms on the right hand side vanish
upon projection to $H^*_c(N_R(K))$, so $[\xi]\in Ker(H^*_c(N_{R+2D}(K))
\ra H^*_c(N_R(K))$.

We now check that $P_\bullet$ is an $\om$-approximate epimorphism.
Let 
$$
[\si]\in Im(H_{n-*}(X,\ol{X-N_{R+2D}(K)})\ra H_{n-*}(X,\ol{X-N_{R}(K)})), 
$$
then $\si$ lifts to a chain $\tau\in C_{n-*}(X)$
so that $\D\tau\in C_{n-*}(\ol{X-N_{R+2D}(K)})$. Let $[\tau]\in H_{n-*}(X, Y_{R+2D})$ 
be the corresponding relative 
homology class. 
Applying $P$ and the chain homotopy $\bar\Phi$, we get
$$
P(\bar P(\tau))-\tau =\D\bar\Phi(\tau)+\bar\Phi(\D \tau).
$$
Since $\bar\Phi(\D\tau)$ vanishes in $C_{n-*}(X,\ol{X-N_R(K)})$,
we get that
$$
[\si]=P_{R+D}(\bar P_{R+2D}([\tau])).
$$
The proof of the last assertion about $\{H^n_c(N_R(K))\}_{R\ge 0}$ and 
$\{H_n(Y_R)\}_{R\ge 0}$  follows since they are approximately isomorphic to zero systems 
$H_0(X, Y_R)$ and $H^0(X, N_R(K))$. 
 \qed

\subsection{Coarse Alexander duality and coarse Jordan separation}
\label{cad}

In this section as in the previous one, we extend complexes indexed by the
nonnegative integers to complexes indexed by $\Z$, by setting the 
remaining groups equal to zero.

 Let $X$, $K$, $D$, $Y_R$, and $\om$ be as in the preceding section.
Composing the morphisms $P_\bul$ and $\bar P_\bul$ with the 
boundary operators for long exact sequences of pairs, we
obtain the compositions $A_{R+D}$
\begin{equation}
\label{f1}
H^*_c(N_{R+D}(K))\stackrel{P_{R+D}}{\RA} H_{n-*}(X,Y_R)
\stackrel{\D}{\simeq}\tilde H_{n-*-1}(Y_R)
\end{equation}
and $\bar A_{R+D}$
\begin{equation}
\label{f2}
\tilde H_{n-*-1}(Y_{R+D})
\stackrel{\D^{-1}}{\simeq}H_{n-*}(X,Y_{R+D})
\stackrel{\bar P_{R+D}}{\RA}H^*_c(N_R(K)).
\end{equation}
Similarly, composing the maps from (\ref{e1})-(\ref{e2}) with
boundary operators and their inverses, we get:
\begin{equation}
\label{f3}
H^*_c(Y_{R})\stackrel{A_{R}}{\lra} \tilde H_{n-*-1}(N_{R+D}(K))
\end{equation}
and
\begin{equation}
\label{f4}
\tilde H_{n-*-1}(N_{R}(K))
\stackrel{\bar A_{R}}{\lra}H^*_c(Y_{R+D}).
\end{equation}

\begin{theorem}[Coarse Alexander duality]
\label{aduality}

1. The morphisms $A_\bul$ and $\bar A_\bul$ in 
(\ref{f1})-(\ref{f4}) are $(\om,\om)$-approximate
isomorphisms.

2.   The maps
$A_\bullet$ in (\ref{f1}) and (\ref{f3}) have  displacement at most $D$. 
The map $\bar A_\bullet$ in (\ref{f2}) (respectively (\ref{f4}))
has  displacement at most $D$ in the sense that if
$\si\in Z_{n-*-1}(Y_{R+D})$  ($\si\in Z_{n-*-1}(N_R(K))$, and
$\si=\D\tau$ for $\tau\in C_{n-*}(X)$, then
the support of $\bar A_{R+D}([\si])$  (respectively $\bar A_R([\si])$)
is contained
in $N_D(Support(\tau))$.   
\end{theorem}

\no
Like ordinary Alexander duality, this theorem follows directly 
from Theorem \ref{Coarse Poincare duality}, 
and the long exact sequence for pairs.

Combining Theorem \ref{aduality} with Corollary \ref{retractionres}
we obtain:

\begin{theorem}
[Coarse Alexander duality for $FP_k$ groups]
\label{adualityres}
Let $X$ be a coarse $PD(n)$ space, and let $G$, $P_*$, $G\acts X$,
$f$, and $K$ be as in the statement of Corollary \ref{retractionres}.
Then

1.  The  family of compositions 
$$
\tilde H_{n-i-1}(Y_{R+D})
\stackrel{\bar A}{\ra}H_c^i(N_R(K))
\stackrel{f_R^i}{\lra} H^i(P_*;\Z G)
$$
defines an approximate isomorphism 
when $ i<k$, and an approximate monomorphism 
when $i=k$.  Recall that for $i<k$ we have
a natural isomorphism $H^i(P_*,\Z G)\simeq H^i(G,\Z G)$.

2. The family of compositions 
$$
\tilde H_{i} (P_*; \Z G) \to \tilde H_{i}(N_R(K)) 
\stackrel{\bar A_R}{\lra} H^{n-i-1}_c (Y_{R+D})
$$
is an approximate isomorphism 
when $ i<k$, and an approximate epimorphism when $i=k$.  
Recall that $\t H_i(P_*;\Z G)=\{0\}$ for $i<k$ since $G$  is
of type $FP_k$.

\end{theorem}

\begin{theorem}[Coarse Alexander duality for maps]
\label{adualitymapping}
\sloppy{Suppose $X$ is a coarse $PD(n)$  space, $X'$ is
a bounded geometry uniformly $(k-1)$-acyclic metric simplicial complex, and
$f:C_*(X')\ra C_*(X)$ is a uniformly proper chain map.  Let 
$K:= Support(f(C_*(X'))$, $Y_R:= \ol{X- N_R(K)}$. 
Then: }

1.  The  family of compositions 
$$
\tilde H_{n-i-1}(Y_{R+D})
\stackrel{\bar A}{\ra}H_c^i(N_R(K)) \stackrel{H^i_c(f_R)}{\RA} H_c^i(X')$$
defines an approximate isomorphism 
when $i<k$, and an approximate monomorphism 
when $i=k$.  

2. The family of compositions 
$$
\tilde H_{i} (X') \to \tilde H_{i}(N_R(K)) \stackrel{\bar A_R}{\lra} 
H^{n-i-1}_c (Y_{R+D})
$$
is an approximate isomorphism 
when $i<k$, and an approximate epimorphism when $i=k$.\footnote{  
The function $\om$ for the above approximate isomorphisms depends only 
on the distortion of $f$, the acyclicity functions for $X$ and $X'$, and the bounds
on the geometry of $X$ and $X'$.} 

3. Furthermore, these approximate isomorphisms approximately
respect support in the following sense.  There is a function $\zeta:\N\ra\N$
so that if  $i<k$, $S\subset X'$ is subcomplex, $T\defeq Support(f_*(C_*(S)))\subset X$
is the corresponding subcomplex of $X$, and
$\al\in Im(H^i_c(X',\ol{X'-S})\ra H^i_c(X'))$, 
then $\al$ belongs to the image of the composition
$$
\tilde H_{n-i-1}(Y_R\cap N_{\zeta(R)}(T))
\ra \tilde H_{n-i-1}(Y_R)\stackrel{H^i_c(f)\circ \bar A}{\RA}H^i_c(X').
$$  
4. If $k=n+1$, then $H^n_c(X')=\{0\}$ unless $N_R(K)=X$ for some $R$.   
\end{theorem}
\proof Parts 1, 2 and 3 of Theorem follow from Lemma \ref{retractionmap} and Theorem \ref{aduality}. 
Part 4 follows since for $i=n$,  $\{\tilde H_{n-i-1}(Y_{R+D})\}=\{0\}$ is approximately isomorphic 
to the constant system $\{H^n_c(X')\}$.  \qed 

We now give a number of corollaries of coarse Alexander duality.

\begin{corollary}[Coarse Jordan separation for maps]
\label{cjs}
\sloppy{Let $X$ and $X'$ be $n$-dimensional and $(n-1)$-dimensional 
coarse Poincar\'e duality spaces respectively, and let $g:X'\to X$ be a
uniformly proper simplicial map. Then}

1. $g(X')$ coarsely separates $X$ into (exactly) two components.  
 
2. For every $R$, each point of 
 $N_R(g(X'))$ lies within uniform distance from each of the
 deep components of $Y_R:= \ol{X-N_R(g(X'))}$.
 
3.  If $Z\subset X'$ and $X'\not\subset N_R(Z)$ for any $R$, then
$g(Z)$ does not coarsely separate $X$.
\end{corollary}
\proof  We have the following diagram:
$$
\begin{array}{ccc}
\t H_0 (Y_{R}) & \stackrel{H^{n-1}_c(g) \circ \bar A }{\RA} H^{n-1}_c(X')= \Z\\
\uparrow & ~ \\ 
\underset{\underset{R}{\longleftarrow}}{\lim} ~ \t H_0 ^{Deep}(Y_{R})
\end{array}
$$
where the family of morphisms $H^{n-1}_c(g) \circ \bar A$ gives rise to an approiximate isomorphism. 
Thus 
$$
\underset{\underset{R}{\longleftarrow}}{\lim} ~ \t H_0 ^{Deep}(Y_{R})= \Z
$$
which implies 1. Let $x\in N_R(K)$. Then there exists a representative $\al$ 
of a generator of $H^{n-1}_c(X')$ such that $H^{n-1}_c(g)(\al)\in C^{n-1}_c(X)$ is supported uniformly close to 
$x$. We apply Part 3 of Theorem \ref{adualitymapping} to the class $[H^{n-1}_c(g)(\al)]$  to prove 2.  
\qed

As a special case of the above corollary we have: 

\begin{corollary}
[Coarse Jordan separation for submanifolds]
Let $X$ and $X'$ be $n$-dimensional and $(n-1)$-dimensional 
uniformly acyclic PL-manifolds respectively, and let $g:X'\to X$ be a
uniformly proper simplicial map. Then the assertions 
1, 2 and 3 from the preceeding theorem hold. 
\end{corollary}

Similarly to the Corollary \ref{cjs} we get: 

\begin{corollary}[Coarse Jordan separation for groups]
\label{cjsgps}
\sloppy{Let $X$ be a  coarse $PD(n)$-space and $G$ be a $PD(n-1)$-group acting freely simplicially 
on $X$. Let $K\subset X$ be a $G$-invariant subcomplex with $K/G$ compact. 
Then:}

1. $G$ coarsely separates $X$ into (exactly) two components.  
 
2. For every $R$, each point of 
 $N_R(K)$ lies within uniform distance from each of the
 deep components of $\ol{X-N_R(K)}$.
 \end{corollary}

\begin{lemma}
\label{uniquecyclic}
Let $W$ be a bounded geometry metric simplicial complex which is homeomorphic
to a union of $W=\cup_{i\in I} W_i$ of $k$ half-spaces $W_i\simeq\R^{n-1}_+$ 
along their
boundaries.  Assume that for $i\neq j$, the union $W_i\cup W_j$
is uniformly acyclic and is uniformly properly embedded in
$W$.  Let $g:W\ra X$ be a uniformly proper map of $W$ into
a coarse $PD(n)$ space $X$.  Then $g(W)$ coarsely separates
$X$ into $k$ components.  Moreover, there is a unique cyclic
ordering  on the index set $I$ so that for $R$ sufficiently large,
the frontier of each deep component $C$ of $X-N_R(g(W))$
is at finite Hausdorff distance from $g(W_i)\cup g(W_j)$
where $i$ and $j$ are adjacent with 
respect to the cyclic ordering.
\end{lemma}
\proof
We have $H^{n-1}_c(W)\simeq \Z^{k-1}$, so, arguing analogously to Corollary \ref{cjs}, we see that 
$g(W)$ coarsely separates $X$ into $k$ components.  Applying coarse
Jordan separation and the fact that no $W_i$ coarsely
separates $W_j$ in $W$, we can define the desired
cyclic ordering by declaring that $i$ and $j$
are consecutive iff $g(W_i)\cup g(W_j)$ coarsely
separates $X$ into two deep components (Corollary \ref{cjs}), one of
which is a deep component of $X-g(W)$.  We leave the
details to the reader.
\qed

\begin{lemma}
\label{compsstable}
 Suppose $G$ is a group of type $FP_{n-1}$ of cohomological
dimension $\leq n-1$, and let $P_*$, $f$,  $G\acts X$, $K\subset X$
and $Y_R$ 
be as in Theorem \ref{adualityres}.  Then 
every deep component of $Y_R$
is stable for $R\geq D$; in particular, there are
only finitely many deep components of $Y_R$
modulo $G$. If $dim(G)<n-1$ then there is only
one deep component.
\end{lemma}
\proof
The composition
\begin{equation}
\label{invlim5}
\underset{\underset{R}{\longleftarrow}}{\lim}\ \t H^{Deep}_0(Y_R)
\ra \t H^{Deep}_0(Y_D)\stackrel{f^i_D\circ\bar A_D}{\RA} H^{n-1}(P_*;\Z G)
\end{equation} 
is an isomorphism by Theorem \ref{adualityres}.  Therefore
$$
\t H^{Deep}_0(Y_R)
\ra \t H^{Deep}_0(Y_D)
$$
is a monomorphism for any $R\geq D$, and hence every deep component of
$Y_D$ is stable.   If $dim(G)<n-1$ then $H^{n-1}(P_*,\Z G)=\{0\}$, and 
by (\ref{invlim5}) we conclude that $Y_D$ contains only one deep component.
\qed

Another consequence of coarse Jordan separation
is:

\begin{corollary}
\label{pdn-1peripheral}
Let $G\acts X$ be a free simplicial action of a group $G$
of type $FP$ on a coarse $PD(n)$ space $X$,
 and let $K\subset X$ be a $G$-invariant
subcomplex on which $G$ acts cocompactly.  By Lemma \ref{compsstable}
there is an $R_0$ so that all deep components of $X-N_{R_0}(K)$
are stable; hence we have a well-defined collection of 
deep complementary components $\{C_\al\}$ and their  stabilizers $\{H_\al\}$.
If $H\subset G$ is
a $PD(n-1)$ subgroup, then one of the following holds:

1.  $H$ coarsely separates $G$.

2.  $H$ has finite index in $G$, and so $G$ is a $PD(n-1)$ group.

3.  $H$ has finite index in $H_\al$ for some $\al$.

In particular, $G$ contains only finitely many conjugacy classes of maximal,
coarsely nonseparating $PD(n-1)$ subgroups.
\end{corollary}
\proof
We assume that $H$ does not coarsely separate $G$.  Pick a basepoint 
$\star\in K$, and let $W\defeq H(\star)$ be the $H$-orbit of $\star$.
Then by Corollary \ref{cjsgps} there is an $R_1$ so that 
$X-N_{R_1}(W)$ has two deep components $C_+,\,C_-$ and both are stable.
Since $H$ does not coarsely separate $G$, we may assume that
  $K\subset N_{R_2}(C_-)$ for some $R_2$.   Therefore $C_+$
has finite Hausdorff distance from some deep component $C_\al$ of
$X-N_{R_0}(K)$, and clearly the Hausdorff distance between the
frontiers $\D C_+$ and $\D C_\al$ is finite.  Either $H$ preserves $C_+$ and $C_-$, or it
contains an element $h$ which exchanges the two.  In the latter
case, $h(C_\al)$ is within finite Hausdorff distance from
$C_-$; so in this case $K$ is contained in $N_r(W)$ for some $r$,
and this implies 2.   When $H$ preserves $C_+$ then we have
$H\subset H_\al$, and since $H$ acts cocompactly on $\D C_+$,
it also acts cocompactly on $\D C_\al$ and hence $[H_\al:H]<\infty$.
\qed

\subsection{The proofs of Theorems \ref{mainduality} and \ref{main3d}}
\label{mainargument}

{\bf Sketch of the proof of Theorem \ref{mainduality}.}  Consider an action $G\acts X$
as in the statement of Theorem \ref{mainduality}.  Let
$K\subset X$ be a $G$-invariant subcomplex with $K/G$ compact.
By Lemma \ref{compsstable} the deep  components
of $X-N_R(K)$ stabilize at some $R_0$, and hence we have
a collection of deep  components
$C_\al$ and their  stabilizers
$H_\al$.  Naively one might hope that for some $R\geq R_0$,
the tubular neighborhood
$N_R(K)$ is acyclic, and
the frontier of $N_R(K)$ breaks up into connected components
which are in one-to-one correspondence with the $C_\al$'s,
each of which is  acyclic and has the same compactly supported
cohomology as $\R^{n-1}$.  Of course, this is too much to hope for,
but there is a coarse analog which does hold.  To explain this
we first note that the systems $\t H_*(N_R(K))$ and  $H^*_c(N_R(K))$ are 
approximately zero and approximately constant respectively
by Corollary \ref{retractionres}.
Applying coarse Alexander duality, we find that the
systems $H^*_c(Y_R)$ and $\t H_*(Y_R)$
corresponding to the complements $Y_R\defeq \ol{X-N_R(K)}$ 
are approximately zero
and approximately constant, respectively.  Instead of 
looking at the frontiers of the neighborhoods $N_R(K)$,
we look at metric annuli $A(r,R)\defeq \ol{N_R(K)-N_r(K)}$
for $r\leq R$.  One can try to compute the (co)homology
of these annuli using a Mayer-Vietoris sequence for the
covering $X=N_R(K)\cup Y_r$; however, the input to this
calculation is only approximate, and the system of annuli
does not form a direct or inverse system in any useful way.
Nonetheless, there are finite direct systems of nested annuli
of arbitrary depth for which one can understand the
(co)homology, and this allows us\footnote{There is an extra
complication in calculating $H^{n-1}_c$ for the annuli
which we've omitting from this sketch.} to apply results from 
section \ref{recognitionsection} to see that the $H_\al$'s
are Poincare duality groups.

\bigskip
\no
{\bf The proof of Theorem \ref{mainduality}.}
We now assume that $G$ is a group of type $FP$ acting
freely simplicially on a coarse $PD(n)$ space $X$.  
This implies that $dim(G)\leq n$, so by Lemma \ref{shortres} 
there is a resolution
$0\ra P_n\ra \ldots \ra P_0\ra \Z\ra 0$ 
 of $\Z$ by finitely generated projective
$\Z G$-modules.  We may construct
$G$-equivariant (augmentation preserving) 
chain mappings $\rho:C_*(X)\ra P_*$
and $f:P_*\ra C_*(X)$ using the acyclicity of $C_*(X)$ and
$P_*$; the composition $\rho\circ f:P_*\ra P_*$
is $\Z G$-chain homotopic to the identity.   If  $L\subset X$ is a 
$G$-invariant subcomplex for which $L/G$ is compact,
then we get an induced homomorphism 
$$
H^*(G;\Z G)\stackrel{H^*(\rho)}{\lra}H^*(X/G;\Z G)\ra H^*(L/G;\Z G)\simeq 
H^*_c(L);$$
 abusing notation we will denote this composition
by $H^*(\rho)$.   

Let $K\subset X$ be
a connected, $G$-invariant subcomplex so that $K/G$
is  compact and the image of $f$ is supported in $K$.
For $R\geq 0$ set $Y_R\defeq \ol{X-N_R(K)}$.
Corollary \ref{retractionres} tells us that the families of maps\begin{equation}
\label{1}
\{0\}\ra \{\t H_*(P_*;\Z G)\}\ra \{\t H_*(N_R(K))\}
\end{equation}
 \begin{equation}
\label{2}
H^*_c(f):H^*_c(N_R(K))\ra H^*(G;\Z G)\simeq H^*(P;\Z G).
\end{equation}
define approximate isomorphisms.
Applying 
Theorems   \ref{adualityres} 
 we get approximate isomorphisms
 \begin{equation}
\label{h*clemma}
\mbox{$\{0\}\ra H^k_c(Y_R)$ for all $k$}
\end{equation}
 and
 \begin{equation}
\label{h*lemma}
\mbox{$\phi_{k,R}:\t H_k(Y_{R})\ra H^{n-k-1}(P_*;\Z G)\simeq H^{n-k-1}(G;\Z G)$
for all $k$}.
\end{equation}
We denote $\phi_{*,D}$ by $\phi_*$. 

We now apply Lemma \ref{compsstable} to see that every 
deep component of $X-N_D(K)$ is stable.  Let $\{C_\al\}$
denote the collection of deep components of 
$X-N_D(K)$, and set $Y_{R,\al}\defeq Y_R\cap C_\al$
and $Z_{R,\al}\defeq \ol{X-Y_{R,\al}}$.   Note that 
for every $\al$, and $D<r<R$ we have
$Z_{R,\al}\cap Y_{r,\al}=\ol{N_{R}(K)-N_{r}(K)}\cap C_\al$.

\begin{lemma}
\label{h0lemma}
1. There is an $R_0$ so that if $R\geq R_0$ then
$Y_{R,\al}=\ol{X-Z_{R,\al}}$ and $Z_{R,\al}=N_{R-R_0}(Z_{R_0,\al})$.

2. The systems $\{\t H_k(Y_{R,\al})\}$, $\{\t H_k(Z_{R,\al})\}$,
$\{H^k_c(Y_{R,\al})\}$, $\{H^k_c(Z_{R,\al})\}$ are approximately
zero for all $k$.
\end{lemma}
\proof
Pick $R_0$ large enough that all shallow components of
$X-N_D(K)$ are contained in $N_{R_0-1}(K)$.  Then for all $R\geq R_0$,
$\D C_\al\cap Y_R=\emptyset$ and hence $Y_{R,\al}$,
like $Y_R$ itself,  is the closure of its interior; this implies
that $Y_{R,\al}=\ol{X-\ol{X-Y_{R,\al}}}=\ol{X-Z_{R,\al}}$.
We also have $Z_{R,\al}=N_R(K)\sqcup (\sqcup_{\be\neq\al} C_\be)$ for
all $R\geq R_0$.  Since $\sqcup_{\be\neq\al}N_R(C_\be)\subset
N_{R_0+R}(K)\cup(\sqcup_{\be\neq\al}C_\be)$,
we get
$$
N_R(Z_{R_0,\al})=N_{R_0+R}(K)\cup(\sqcup_{\be\neq\al} N_R(C_\be))
$$
$$
=N_{R_0+R}(K)\cup(\sqcup_{\be\neq\al}C_\be)
$$
$$
=Z_{R_0+R,\al}.
$$
Thus we have proven 1.

To prove 2, we first note that $\{\t H_0(Y_{R,\al})\}$ is approximately
zero by the stability of the deep components $C_\al$.    When $R\geq R_0$
then $Z_{R,\al}$ is connected (since $N_R(K)$ and each $C_\be$ are connected),  
and this says that  $\{\t H_0(Z_{R,\al})\}$ is approximately zero.   
When $R\geq R_0$ then $Y_R$ is the disjoint union $\sqcup_\al Y_{R,\al}$, 
so we have direct sum decompositions $H_k(Y_R)=\oplus_\al H_k(Y_{R,\al})$
and $H^k_c(Y_R)=\oplus_\al H^k_c(Y_{R,\al})$ which are compatible
projection homomorphisms.  This together with (\ref{h*clemma})
and (\ref{h*lemma}) implies that $\{\t H_k(Y_{R,\al})\}$
and $\{ H^k_c(Y_{R,\al})\}$ are approximately zero for all $k$.
By part 1 and Theorem \ref{aduality} we get that $\{H^k_c(Z_{R,\al})\}$
and $\{\t H_k(Z_{R,\al})\}$ are approximately zero for all $k$.
\qed

\medskip
\begin{lemma}
\label{nestedannuli}
There is an $R_{min}>D$ so that
for any $R\geq R_{min}$ and any integer $M$, there is a sequence
$R\leq R_1\leq R_2\leq ...\leq R_M$ with the following
property.  Let $A(i,j)\defeq \ol{N_{R_j}(K)-N_{R_i}(K)}\subset Y_{R_i}$, 
and $A_\al(i,j)\defeq A(i,j)\cap C_\al$.  Then for each $1< i<j< M$,

1.  The image of
$\tilde H_k(A(i,j))\ra \tilde H_k(A(i-1,j+1))$
maps isomorphically onto $H^{n-k-1}(G;\Z G)$ under the 
composition $\tilde H_k(A(i-1,j+1))\ra
\tilde H_k(Y_D)\stackrel{\phi_k}{\ra}H^{n-k-1}(G;\Z G)$ 
for $0\leq k\leq n-1$.  The homomorphism 
$\tilde H_n(A(i,j))\ra \tilde H_n(A(i-1,j+1))$ is zero.

2.   $H^k(\rho):H^k(G;\Z G)\ra H^k_c(A(i,j))$ maps $H^k(G;\Z G)$ isomorphically
onto the image of $H^k_c(A(i-1,j+1))\ra H^k_c(A(i,j))$
 for $0\leq k<n-1$. 
 
3. There is a system of homomorphisms 
$H^{n-1}_c(A_\al(i,j))\stackrel{\th_{i,j}^\al}{\lra}\Z$ 
(compatible with the 
inclusions $A_\al(i,j)\ra A_\al(i-1,j+1)$)
so that the image of 
$H^{n-1}_c(A_\al(i-1,j+1))\ra H^{n-1}_c(A_\al(i,j))$
maps isomorphically to $\Z$ under $\th_{i,j}^\al$.

4. For each $\al$, $\t H_0(A_\al(i,j))\stackrel{0}{\ra}\t H_0(A_\al(i-1,j+1)).$
\end{lemma}

\proof
We choose $R_{min}$ large enough that for any $R\geq R_{min}$,
the following inductive construction is valid. 
Let $R_1\defeq R$.  Using the approximate isomorphisms (\ref{1}), 
(\ref{2}), (\ref{h*clemma}), (\ref{h*lemma}), and Lemma 
\ref{h0lemma}, we inductively  choose $R_{i+1}$
so that:

A.  $\tilde H_k(N_{R_i}(K))\stackrel{0}{\ra} \tilde H_k(N_{R_{i+1}}(K))$
for $0\leq k\leq n$.

B. $Im(\tilde H_k(Y_{R_{i+1}})\ra \tilde H_k(Y_{R_i}))$
maps isomorphically to $H^{n-k-1}(G;\Z G)$ under $\phi_{k,R_i}$ for
$0\leq k< n$, and $Im(\tilde H_k(Y_{R_{i+1}})\ra \tilde H_k(Y_{R_i}))$
is zero when $k=n$.

C. $Im(H^*_c(N_{R_{i+1}}(K))\ra H^*_c(N_{R_i}(K)))$ maps
isomorphically onto $H^*(G;\Z G)$ under $H^*_c(f)$.

D.  $H^*_c(Y_{R_i})\stackrel{0}{\ra} H^*_c(Y_{R_{i+1}})$.

E.  For each $\al$, $H^{n-1}_c(Y_{R_i,\al})\stackrel{0}{\ra} H^{n-1}_c(Y_{R_{i+1},\al})$,
and $H^{n-1}_c(Z_{R_{i+1},\al})\stackrel{0}{\ra} H^{n-1}_c(Z_{R_{i},\al})$.

F.  For each $\al$,  $\t H_0(Y_{R_{i+1},\al})\stackrel{0}{\ra}\t H_0(Y_{R_i,\al})$ and
$\t H_0(Z_{R_i,\al})\stackrel{0}{\ra}\t H_0(Z_{R_{i+1},\al})$.

Now take $1<i<j<M$, and consider the map of  Mayer-Vietoris sequences for the
decompositions $X=N_{R_j}(K)\cup Y_{R_i}$ and $X=N_{R_{j+1}}(K)\cup
Y_{R_{i-1}}$:
$$\begin{array}{ccccc}
\t H_{k+1} (X) \ra & \t H_k(A(i,j))\ra & \t H_k(N_{R_j}(K))\oplus\t H_k(Y_{R_i}) & \ra & \t H_k(X) \\
\downarrow  & \downarrow &   0\downarrow ~~~~~~~~ \downarrow &   &   \downarrow \\
 \t H_{k+1} (X) \ra & \t H_k(A(i-1,j+1))\ra & \t H_k(N_{R_{j+1}}(K))\oplus\t H_k(Y_{R_{i-1}}) & \ra & \t H_k(X) \\
~  & \downarrow \phi_k\restr_{A(i-1,j+1)}       & ~~~~~~~~~~~~~~~~~~\downarrow \phi_k &       & ~  \\
~ & H^{n-k-1}(G, \Z G) \ra & ~~~~~~~~~~ H^{n-k-1}(G, \Z G) & ~  
\end{array}$$

Since $\t H_*(X)=\{0\}$, conditions A and B and the
diagram imply the first part of assertion 1.   The same
Mayer-Vietoris diagram for $k=n$ implies the second part.

Let $0\leq k< n-1$. Consider the commutative diagram of Mayer-Vietoris sequences:

$$
\begin{array}{cccc}
~ & H^{k}(G, \Z G) \to & H^{k}(G, \Z G) & ~ \\
~  & H^k(\rho)\downarrow   ~~~~~~~~~~~  & H^k(\rho)\downarrow & ~  \\
H^{k}_c (X) \ra & H^{k}_c (N_{R_{j+1}}(K)) \oplus  
H^{k}_c (Y_{R_{i-1}}) \ra & H^{k}_c (A(i-1, j+1)) & 
\ra H^{k+1}_c(X) \\
\downarrow  & \downarrow ~~~~~~~~  0 \downarrow &  \downarrow &   \downarrow \\
 H^{k}_c (X) \ra & H^{k}_c (N_{R_{j}}(K)) 
 \oplus H^{k}_c (Y_{R_{i}}) \ra & H^{k}_c (A(i, j)) & \ra H^{k+1}_c
(X) \\
 
\end{array}$$
Assertion 2 now 
follows from the fact that $H^{k}_c (X)\cong H^{k+1}_c (X)=0$,
conditions C and D, and the diagram.   

Assertion 3 follows from condition E, the fact that
$H^n_c(X)\simeq \Z$, and the following
commutative diagram of Mayer-Vietoris sequences ($\th_{i,j}^\al$
is the coboundary operator in the sequence):

$$
\begin{array}{ccc}
H^{n-1}_c (Z_{R_{j+1},\al}) \oplus  
H^{n-1}_c (Y_{R_{i-1},\al}) \ra & H^{n-1}_c (A_\al(i-1, j+1)) & 
\stackrel{\th_{i-1,j+1}}{\RA} H^{n}_c(X) \to 0\\
0\downarrow ~~~~~~~~~~~ 0\downarrow & \downarrow &  \downarrow  ~~~\\
H^{n-1}_c (Z_{R_{j},\al}) 
 \oplus H^{n-1}_c (Y_{R_{i},\al}) \ra & H^{n-1}_c (A_\al(i, j)) & 
 \stackrel{\th_{i,j}}{\RA} H^{n}_c(X) \to 0\\ 
\end{array}$$

Assertion 4 follows from condition F and the following commutative
diagram:

$$\begin{array}{ccccc}
\t H_{1} (X) \ra & \t H_0(A_\al(i,j))\ra & \t H_0(Z_{R_j,\al})\oplus\t H_0(Y_{R_i,\al}) & \ra & \t H_0(X) \\
\downarrow  & \downarrow &   0\downarrow ~~~~~~~~ 0\downarrow &   &   \downarrow \\
 \t H_{1} (X) \ra & \t H_0(A_\al(i-1,j+1))\ra & \t H_0(Z_{R_{j+1},\al})\oplus\t H_0(Y_{R_{i-1},\al}) & \ra & \t H_0(X) \\
\end{array}$$
\qed

\begin{corollary}
If $G$ is an $(n-1)$-dimensional duality group, then each deep component stabilizer
is a $PD(n-1)$ group.
\end{corollary}
\proof Fix a deep  component $C_\al$ of $X-N_D(K)$, and let $H_\al$ be its stabilizer in $G$. 
Let $R=D$, $M=4k+2$, and apply the construction of Lemma
\ref{nestedannuli} to get $D\leq R_1\leq R_2\leq\ldots\leq R_{4k+2}$
satisfying the conditions of Lemma \ref{nestedannuli}.  

Pick $1<i<j<M$.  The mappings 
$\t H_\ell(A(i,j))\to \t H_\ell (A(i-1,j+1))$ are zero for each $\ell =1,...,n$
by part 1 of Lemma \ref{nestedannuli}, since  $H^k(G, \Z G)= 0$ for $k< n-1$.
Because $A(p,q)$ is the disjoint union  $\amalg_\al A_\al(p,q)$ for all $0<p<q<M$,
we actually have $\t H_\ell(A_\al(i,j))\stackrel{0}{\to} \t H_\ell (A_\al(i-1,j+1))$
for $1\leq \ell\leq n$.  By part 4 of Lemma \ref{nestedannuli} the same
assertion holds for $\ell =0$.  Applying 
Theorem \ref{finitenesstheorem} to the
chain complexes $C_*(A_\al(i,j))$, we see that 
when $k>2n+5$, $H_\al$ is a group
of type $FP(n)$.  Since $dim(H_\al)\leq dim(G)=n-1$ it follows that 
$H_\al$ is of type $FP$ (see section  \ref{groupprelim}).

The mappings $H^\ell _c(A_\al(i-1,j+1))\to H^\ell _c(A_\al(i,j))$
are zero for $0\leq \ell<n-1$ by part 2 of Lemma \ref{nestedannuli}
and the fact that $A(p,q)=\amalg_\al A_\al(p,q)$.  By parts 1 and
2 of Lemma \ref{finlem}, we have $H^k(H_\al,\Z H_\al)=\{0\}$
for $0\leq k<n-1$, and $H^{n-1}(H_\al,\Z H_\al)\simeq\Z$ by part
3 of Lemma \ref{nestedannuli}.  Hence $H_\al$ is a $PD(n-1)$ group.
\qed

\no
{\bf Remark.}  For the remainder of the proof, we really only need
to know that each deep component stabilizer is of type $FP$.  

{\em Proof of Theorem \ref{mainduality} concluded.}
Let $C_1,\ldots,C_k$ be  a set of representatives
for the $G$-orbits of deep  components of
$X-N_R(K)$, and let $H_1,\ldots,H_k\subset G$ denote their
stabilizers.   Since $G$ and each $H_i$ is of type $FP$, the
group pair $(G,\{H_i\})$ has finite type (section \ref{groupprelim}).
By Lemma \ref{relretractionres}, we have 
$$
H^*(G,\{H_i\};\Z G)\simeq \underset{\underset{R}{\lra}}{\lim}~ H^*_c(X,Y_R),
$$
while $\lim_R H^*_c(X,Y_R)\simeq \lim_R  H_{n-*}(N_R(K))$ by 
Coarse Poincare duality,  and 
$$
\underset{\underset{R}{\lra}}{\lim}~ H_*(N_R(K)) \simeq H_*(X)\simeq H_*(pt)
$$
 since homology commutes with direct limits.
Therefore the group pair $(G,\{H_i\})$ satisfies one of the
criteria for $PD(n)$ pairs (see section \ref{groupprelim}),
 and we have proven Theorem \ref{mainduality}.
 \qed

\medskip
We record a variant of Theorem \ref{mainduality} which describes the
geometry of the action $G\acts X$ more explicitly:

\begin{theorem}
\label{peripheralversion}
Let $G\acts X$ be as in Theorem \ref{mainduality}, and let 
$K\subset X$ be a $G$-invariant subcomplex with $K/G$ compact.
Then there are $R_0,\,R_1,\,R_2$ so that

1.  The deep components $\{C_\al\}_{\al\in I}$ of $X-N_{R_0}(K)$ are
all stable, there are only 
finitely many of them modulo $G$, and their stabilizers $\{H_\al\}_{\al\in I}$
are $PD(n-1)$ groups.

2.  For all $\al\in I$, the frontier $\D C_\al$ is connected, and 
$N_{R_1}(\D C_\al)$ has precisely two deep complementary components,
$E_\al$ and $F_\al$,
where $E_\al$ has Hausdorff distance at most $R_2$ from $C_\al$.
Unless $G$ is a $PD(n-1)$ group,  the distance function
$d(\D C_\al,\cdot)$ is unbounded on $K\cap F_\al$. 

3.  The Hausdorff distance between $X- \amalg_\al E_\al$ and $K$ is at most 
$R_2$.

\end{theorem}
\proof This is clear from the discussion above.
\qed

\no
We remark that there are $\al_1\neq\al_2\in I$ so that the Hausdorff
distance $$d_H(\D C_{\al_1},\D C_{\al_2})<\infty$$ iff $G$ is a $PD(n-1)$ group.

\begin{lemma}
\label{sim}
Let $G\acts X$ be as Theorem \ref{mainduality}, and let $K,\,C_\al,\,H_\al,C_i,\,H_i$ 
be  as in the conclusion of the proof of Theorem \ref{mainduality}.   
If $X$ is simply connected and the groups $H_i$ admit
finite $K(H_i,1)$'s, then $G$ admits a finite $K(G,1)$.
There exists a contractible coarse $PD(n)$ space $X'$ on which $G$ acts freely and 
simplicially with the following properties: 

1. There is a $G$-equivariant proper homotopy equivalence $\phi:X\to X'$ which is 
a homeomorphism away from a finite tubular neighborhood of $K$.   

2. There is a contractible subcomplex $K'\subset X'$ on which $G$ acts cocompactly. 
All components of $X'- K'$ are deep and stable. 

3. The mapping $\phi$ induces a bijection between the deep components $C_\al$ and 
components of $X'- K'$. 
\end{lemma}
\proof For each $1\leq i\leq k$, let $W_i$ be the universal cover
of a finite Eilenberg-MacLane space for $H_i$, and specify
an  $H_i$-equivariant map  $\psi_i:\D C_i\ra W_i$, where
$\D C_i$ is the frontier of $C_i$. 
  We can $G$-equivariantly
identify  the disjoint union $\amalg_{\al\in G(i)}\D C_\al$ with
the twisted product $G\times_{H_i}\D C_i$, and obtain an
induced $G$-equivariant mapping 
$$
\Psi: \cup_\al \,\D C_\al=\cup_i(\amalg_{\al\in G(i)}\,\D C_\al)
\ra \amalg_i \,(G\times_{H_i}W_i).$$
Let $K^+:= \ol{X-\amalg_\al C_\al}$. We now cut $X$ open
along the disjoint union $\D C\defeq \amalg_\al \,\D C_\al$ to obtain
a new complex 
$$
\ddot X:= K^+ \amalg (\amalg_\al C_\al)
$$ 
which contains two copies $\D_+ C\subset K^+$ and $\D_-C\subset \amalg_\al C_\al$
of $\D C$. We let $\Psi_\pm$ be the corresponding copies of the mapping $\Psi$. 
Now define $K'$ as the union (along $\D_+ C$) of $K^+$ and the mapping cylinder of
 $\Psi_+$  and define $Y'$ as the union (along $\D_- C$) of $\amalg_\al C_\al$ and 
the mapping cylinder of $\Psi_-$.   
Finally obtain $X'$ gluing $K'$ and $Y'$ 
along the copies of $W:=\amalg_i (G\times_{H_i}W_i)$.  
The group $G$ still acts on $X'$ freely and simplicially and clearly $K'/G$ is compact. By applying 
Van-Kampen's theorem and Mayer-Vietoris sequences, it follows
that $X'$ and $K'$ are uniformly contractible.  
Assertion 1 is clear from the construction of $X'$.  The remaining
assertions follow easily from the first.
\qed

\medskip
\no
{\bf Proof of Theorem \ref{main3d}.} By the main theorem deep components stabilizers $H_i$ 
are $PD(2)$-groups and hence are surface groups by 
\cite{EM,EL}. 
Theorem \ref{main3d} now follows by applying Lemma \ref{sim} (where the complexes $W_i$ in 
the proof are homeomorphic to $\R^2$). 
\qed

In  Proposition \ref{peripheraluniqueness} 
 we generalize the uniqueness theorem of the peripheral  
structure from  3-dimensional  manifolds  to  $PD(n)$  pairs.  

\begin{theorem}\label{jsj}
(Jaco and Shalen \cite{jacoshalen}, Johannson  \cite{johannson}.) 
Let $M$ be a compact connected acylindrical 3-manifold with aspherical incompressible 
boundary components $S_1,\ldots,S_m$.  Let $N$ be a compact 3-manifold 
homotopy-equivalent to $M$, with incompressible boundary
components $Q_1,\ldots,Q_n$, and $\varphi: \pi_1(M)\to \pi_1(N)$ 
be an isomorphism. Then $\varphi$ preserves the peripheral structures of $\pi_1(M)$ and $\pi_1(N)$ in the 
following sense. There is a bijection $\beta$ between the set of boundary components of $M$ and the set of 
boundary components on $N$ so that after relabelling via  $\beta$ we have:

$\varphi(\pi_1(S_i))$ is conjugate to $\pi_1(Q_i))$ in $\pi_1(N)$. 
\end{theorem}

\begin{proposition}
\label{peripheraluniqueness}
Let $(G,\{H_i\}_{i\in I})$ be a $PD(n)$ pair, 
where $G$ is not a $PD(n-1)$ group, and 
$H_i$ does not coarsely separate $G$ for any $i$.
Now let $G\acts X$
be a free simplicial action on a coarse $PD(n)$ space,
and let $(G,\{L_j\}_{j\in J})$ be the group pair obtained by applying
Theorem \ref{mainduality} to this action.  
Then  there is a bijection
$\be:I\ra J$ so that $H_i$ is conjugate to $L_{\be(i)}$ for
all $i\in I$.
\end{proposition}
\proof  Under the assumptions above, each $H_i$ and $L_j$
is a maximal $PD(n-1)$ subgroup (see Lemma \ref{hi'sdistinct}).
By Corollary \ref{pdn-1peripheral}, each $H_i$ is conjugate to
some $L_j$, and by Lemma~\ref{hi'sdistinct} this defines an
injection $\beta:I\ra J$.  Consider the double $\hat G$ of $G$
over the $L_j$'s.  Then the double of $G$ over the $H_i$'s
sits in $\hat G$, and the index will be infinite unless
$\beta$ is a bijection.  
 \qed

We  now  establish  a  relation  between  the  acylindricity   
assumption  in Theorem  \ref{jsj}  and  coarse  nonseparation  assumption  
in Proposition  \ref{peripheraluniqueness}.  We  first  note  that   if  $M$  
is a  compact  3-manifold  with incompressible  aspherical  
boundary components $S_1,\ldots,S_m$,  then  $M$  is   acylindrical  iff 
$\pi_1(S_i) \cap g (\pi_1(S_j)) g^{-1}=\{e\}$ whenever 
$i\ne j$ or $i=j$ but $g\notin \pi_1(S_i)$.

\begin{lemma}
\label{a}
Suppose $G$ is a duality group and $G\acts X$
is a free simplicial action on a coarse $PD(n)$ space,
and let $(G,\{H_j\}_{j\in J})$ be the group pair obtained by applying
Theorem \ref{mainduality} to this action. Assume that $H_i \cap (g H_j g^{-1})=\{e\}$ 
whenever $i\ne j$ or $i=j$ but $g\notin H_i$. Then no $H_i$ coarsely separates $G$. 
\end{lemma}
\proof
Let $K_0\subset X$ be a connected $G$-invariant subcomplex so that
$K_0/G$ is compact and all deep components of $X-K_0$ are stable.
Now enlarge $K_0$ to a subcomplex $K\subset X$ by throwing in the
shallow (i.e. non-deep) components of $X-K_0$; then $K$ is connected,
$G$-invariant, $K/G$ is compact, and all components of $X-K$ are deep
and stable.  Let $\{C_\al\}$ denote the components of $X-K$, and
let $C_i$ be a component stabilized by $C_i$.  We will show that 
$\D C_i$ does not coarsely separate $K$ in $X$.  Since $K\hook X$ is a
uniformly proper embedding, $G\acts K$ is cocompact, and
$H_i\acts \D C_i$ is cocompact, this will imply the lemma.

For all components $C_\al$ and all $R$, the intersection $H_i\cap H_\al$
acts cocompactly on $N_R(\D C_i)\cap\bar C_\al$, where
$H_\al$ is the stabilizer of $C_\al$; when $\al\neq i$ the
group $H_i\cap H_\al$ is trivial, so in this case $Diam(N_R(\D C_i)\cap\bar C_\al)<\infty$.
For each $R$ there are only finitely many $\al$ -- modulo $H_i$ --
for which $N_R(\D C_i)\cap C_\al\neq\emptyset$, so there is
a constant $D_1=D_1(R)$ so that if $\al\neq i$ then 
$Diam(N_R(\D C_i)\cap C_\al)<D_1$.  Each $\D C_\al$ is connected and
$1$-ended, so we have an $R_1=R_1(R)$ so that if $\al\neq i$, and 
$x,\,y\in\D C_\al-N_{R_1}(\D C_i)$,
then $x$ may be joined to $y$ by a path in $\D C_\al-N_R(\D C_i)$.

By Corollary \ref{cjsgps}, there is a function $R_2=R_2(R)$ so that
if $x,\,y\in K-N_{R_2}(\D C_i)$ then $x$ may be joined to 
$y$ by a path in $X-N_R(\D C_i)$.  

Pick $R$, and let $R'=R_2(R_1(R))$.   If $x,\,y\in K-N_{R'}(\D C_i)$
then they are joined by a path $\al_{xy}$ in $X-N_{R_1(R)}(\D C_i)$.
For each $\al\neq i$, the portion of $\al_{xy}$ which enters
$C_\al$ may be replaced by a path in $\D C_\al-N_R(\D C_i)$.
So $x$ may be joined to $y$ in $K-N_R(\D C_i)$.  Thus
$\D C_i$ does not coarsely separate $K$ in $X$.
\qed

\begin{lemma}
\sloppy{Let $M$ be a compact 3-manifold with $\D M\neq\emptyset$, with
aspherical 
incompressible nonempty boundary components $S_1,\ldots,S_m$. 
Then $M$ is acylindrical if and only if $\pi_1(M)$ is not a surface group and no $H_i= \pi_1(S_i)\subset \pi_1(M)=G$ 
coarsely separates $G$.  }
\end{lemma}
\proof The implication $\Rightarrow$  follows from Lemma \ref{a}. To establish 
 $\Leftarrow$ assume that $M$ is not acylindrical. This implies that there exists a nontrivial decomposition
of $\pi_1(M)$ as a graph of groups with a single edge group $C$ which is a cyclic subgroup of some $H_i$.  Thus  
 $C$ coarsely separates $G$. Since $[G: H_i]=\infty$ it follows that $H_i$ coarsely separates $G$ as well. \qed

\begin{corollary}
\sloppy{Suppose  $G$ is not a $PD(n-1)$ group, both $(G,\{H_i\}_{i\in I})$ and $(G,\{L_j\}_{j\in J})$ are 
$PD(n)$ pairs, no $H_i$ coarsely separates $G$, 
and each $L_j$ admits a finite Eilenberg-MacLane space. Then  there is a bijection
$\be:I\ra J$ so that $H_i$ is conjugate to $L_{\be(i)}$ for all $i\in I$. 
Thus the peripheral structure of $G$ in this case is unique.} 
\end{corollary}
\proof Under the above assumptions the double $\hat G$ of $G$ with respect 
to the collection of subgroups $\{L_j\}_{ j\in J}$ 
admits a finite Eilenberg-MacLane space $K(\hat G, 1)$. 
Thus we can take as a coarse $PD(n)$-space $X$ the universal cover of  
$K(\hat G, 1)$. Now apply Proposition 
\ref{peripheraluniqueness}. \qed

\subsection{Applications}
\label{appl}

In this section we discuss examples of $(n-1)$-dimensional
groups which cannot
act on coarse $PD(n)$ spaces.

\medskip
\no
{\bf $2$-dimensional groups with positive Euler characteristic.}
Let $G$ be a group of type $FP_2$ with cohomological dimension $2$.
If the $\chi(G)>0$ then $G$ cannot act freely simplicially on
a coarse $PD(3)$ space.  To see this, note that by Mayer-Vietoris some one-ended
free factor $G'$
of $G$  must have $\chi(G')>0$.  If $G'$ acts on a coarse $PD(3)$ space
then $G'$ contains a collection ${\cal H}$ of surface subgroups
so that $(G',{\cal H})$ is a $PD(3)$ pair.  Since the double of a 
$PD(3)$ pair is a $PD(3)$ group ( which has zero Euler characteristic) by
Mayer-Vietoris  we have  $\chi(G')\leq 0$, which is a  contradiction.

\medskip
\no
{\bf Bad products.}
Suppose $G=\prod_{i=1}^k G_i$ where each $G_i$ is a duality group of 
dimension $n_i$, and $G_1,\,G_2$ are not Poincare duality groups.
Then $G$ cannot act freely simplicially on a coarse
$PD(n)$ space, where $n-1=\sum_{i=1}^k n_i$.   

\proof
Let $G\acts X$ be a free simplicial action on a coarse $PD(n)$ space.

{\em Step 1.  $G$ contains a $PD(n-1)$ subgroup.}  This follows
by applying Theorem \ref{mainduality} to $G\acts X$, since otherwise
$G\acts X$ is cocompact and Lemma \ref{cptsuppiso} would give
$H^n(G;\Z G)\simeq \Z$, contradicting $dim(G)=n-1$.

We apply Theorem \ref{mainduality} to see that $G\acts X$ 
defines deep complementary component stabilizers $H_\al\subset G$
which are $PD(n-1)$ groups.

{\em Step 2. Any $PD(n-1)$ subgroup $V\subset G$ virtually splits
as a product $\prod_{i=1}^k V_i$ where $V_i\subset G_i$ 
is a $PD(n_i)$ subgroup.  Consequently each $G_i$ contains a $PD(n_i)$
subgroup.}

\begin{lemma}
A $PD(m)$ subgroup $V$ of a $m$-dimensional product group
$W\defeq \prod_{i=1}^k W_i$ contains a finite index subgroup $V'$
which splits as a product $V'=\prod_{i=1}^k V_i$
where $V_i\subset W_i$ is a Poincare duality group
of dimension $dim(W_i)$.
\end{lemma}

\proof
Look at the kernels of the projections 
$$\hat p_j:W\ra \prod_{i\neq j}W_i$$
restricted to $V$.
The dimension of the middle group in a short
exact sequence 
 has dimension at most the sum of the dimensions
of the other two groups.  Applying this to the exact sequence
$$1\to W_j\cap V\to V\to \hat p_j(V)\to 1$$
we get that $W_j\cap V$ has
the same dimension as $W_j$.
Hence $\prod_j (W_j\cap V)$ has the same dimension
as $V$, so it has finite index in $V$ (see section \ref{groupprelim}).
Therefore $\prod_j(W_j\cap V)$  is a $PD(n)$ group and 
so the factor groups $(W_j\cap V)$ are $PD(dim(W_j))$ groups.
\qed

Step 3.   {\em No $PD(n-1)$ subgroup $V\subset G$ can coarsely separate $G$. }
This follows immediately from step 2 and:

\begin{lemma}
For $i=1,\,2$ let $A_i\subset B_i$ be finitely generated groups, with
$[B_i:A_i]=\infty$.   Then $A_1\times A_2$ does not coarsely separate
$B_1\times B_2$.
\end{lemma}
\proof  Suppose 
that $x=(x_1,x_2), y=(y_1,y_2)$ are points in the Cayley graphs of $B_1, B_2$ 
which are at distance at least $R$ from $A\defeq A_1\times A_2$. 
Without loss of generality we may assume  that 
$d(x_1, A_1)\ge R/2$.  We then pick a point $x_2'\in B_2$ with distance 
at least $R/2$ from $A_2$ and connect $x_2$ to $x_2'$ by a path $x_2(t)$ 
the the Cayley graph of $B_2$.    The  path $(x_1, x_2(t))$ does not 
intersect $N_{\frac{R}{2}}(A)$. Applying similar argument to $y$ we reduce 
the proof to the case where $d(x_i, A_i)\ge R/2$ and 
$d(y_i, A_i)\ge R/2$, $i=1,2$.  Now connect $x_1$ to $y_1$ by a path 
$x_1(t)$, and $y_2$ to $x_2$ by a path $y_2(t)$; 
it is clear that the paths $(x_1(t), x_2), (y_1, y_2(t))$ do not 
intersect $N_{\frac{R}{4}}(A)$. On the other hand, these paths connect $x$ to 
$(y_1, x_2)$ and $y$ to $(y_1, x_2)$. \qed

\medskip
Step 4.  By steps 1 and  2 we know that each $G_i$ 
contains a $PD(n_i)$ subgroup.  Let $L_i\subset G_i$ be a 
$PD(n_i)$ subgroup for $i>1$.    Set 
$L\defeq G_1\times (\prod_{i=2}^k L_i)$.
Observe that $L$ is not a $PD(n-1)$ group since $G_1$
is not a $PD(n_1)$ group.  Therefore no finite index subgroup of
$L$ can be a $PD(n-1)$ subgroup, see section \ref{groupprelim}.

Step 5. Choose a basepoint $\star\in X$.
We now apply Theorem \ref{peripheralversion} to  the action
$L\acts X$ with $K\defeq L(\star)$, and we let $R_i$, $C_\al$, $H_\al$
$E_\al$, and $F_\al$ be as in the Theorem \ref{peripheralversion}.  Since $L$ has
 infinite index in $G$, the distance function
$d(\D C_\al,\cdot)$ is unbounded on $G(\star)\cap E_\al$ for some
$\al\in I$, while part 2 of Theorem \ref{peripheralversion}
implies that   $d(\D C_\al,\cdot)$ is unbounded on $K\cap F_\al$.
Hence $H_\al$ coarsely separates $G$, which contradicts step 3.
\qed

\medskip
\no
{\bf Baumslag-Solitar groups.}  Pick $p\neq\pm q$, and let $G\defeq BS(p,q)$ denote the
Baumslag-Solitar group with the presentation 
\begin{equation}
\label{bspres}
\<a,b\mid ba^pb^{-1}=a^q\>.
\end{equation}
If $G_1$ is a $k$-dimensional duality group  then 
the direct product $G_1\times G$ does not  act freely 
simplicially on a coarse $PD(3+k)$ space.

We will prove this when $G_1=\{e\}$.   The general case can be proved
using straightforward generalization of the argument
given below, once one applies the ``Bad products'' example above to see that
$G_1$ must be a $PD(k)$ group if $G_1\times G$ acts on a coarse
$PD(3+k)$ space.  
Assume that $G\acts X$ is a free simplicial action on a coarse $PD(3)$ space.
Choosing a basepoint $\star\in X$, we obtain a uniformly proper map
$G\ra X$.   

We recall that the  presentation (\ref{bspres}) defines a graph of groups
decomposition of $G$ with one vertex labelled with $\Z$, one oriented
edge labelled with $\Z$, and where the initial and final edge monomorphisms
embed the edge group as subgroups of index $p$ and $q$ respectively.   The Bass-Serre
tree $T$ corresponding to this graph of groups has the following structure.
The action  $G\acts T$ has one vertex orbit and one edge orbit.   For each
vertex $v\in T$, the vertex stabilizer $G_v$ is isomorphic to $\Z$.   
The vertex $v$  has $p$ incoming edges and $q$ outgoing edges;
the incoming (respectively outgoing) edges  are cyclically
permuted by $G_v$ with ineffective kernel the subgroup
of index $p$ (respectively $q$).

Let $\bar \Si$ be the presentation complex corresponding to
the  presentation (\ref{bspres}), and let $\Si$ denote its universal cover.
Then $\Si$ admits a natural $G$-equivariant fibration $\pi:\Si\ra T$,
with fibers homeomorphic to $\R$.   For each vertex $v\in T$, the inverse
image $\pi^{-1}(v)$ has a cell structure isomorphic to the usual cell structure
on $\R$, and $G_v$ acts freely transitively on the vertices.   For each
edge $e\subset T$, the inverse image $\pi^{-1}(e)\subset \Si$ is homeomorphic
to a strip.    The cell structure on the strip may be obtained as follows.
Take the unit square in $\R^2$ with
the left edge subdivided into $p$ segments and the right edge subdivided
into $q$ segments; then glue the top edge to the bottom edge by translation
and take the induced cell structure on the universal cover.   The edge stabilizer
$G_e$ acts simply transitively on the $2$-cells of $\pi^{-1}(e)$.

We may view $\Si$ as a bounded geometry metric simplicial complex by
taking a $G$-invariant triangulation of $\Si$.   Given $k$ distinct
ideal boundary points $\xi_1,\ldots,\xi_k\in \geo T$ and a basepoint $\star \in T$,
we consider the geodesic rays $\ol{\star\xi_i}\subset T$, take the disjoint
union of their inverse images $Y_i\defeq \pi^{-1}(\ol{\star\xi_i})\subset \Si$ and glue
them together along the copies of $\pi^{-1}(\star)\subset \pi^{-1}(\ol{\star\xi_i})$.
The resulting complex $Y$ inherits bounded geometry metric simplicial complex
structure from $\Si$.   The reader will verify the following assertions:

1.  $Y$ is uniformly contractible.

2.  For $i\neq j$, the union $Y_i\cup Y_j\subset Y$ is uniformly contractible and the
inclusion $Y_i\cup Y_j\ra Y$ is uniformly proper.

3.  The natural map $Y\ra\Si$ is uniformly proper.

4.  The cyclic ordering induced on the $Y_i$'s by the uniformly proper
composition $C_*(Y)\ra C_*(\Si)\ra C_*(X)$  (see Lemma \ref{uniquecyclic})
defines a continous $G$-invariant cyclic ordering on $\geo T$.

Let $a$ be the generator of $G_v$ for some $v\in T$.   Setting
$e_k\defeq (pq)^k$, the 
sequence $ g_k\defeq a^{e_k}$ -- viewed as elements in 
$Isom(T)$ -- converges to the identity as $k\ra \infty$.
So the sequence of induced homeomorphisms of the
ideal boundary of $T$ converges to the identity.
The invariance of the cyclic ordering clearly implies that
$g_k$ acts trivially on the ideal boundary of $T$
for large $k$.   This implies that $g_k$ acts
trivially on $T$ for large $k$.   Since
this is absurd, $G$ cannot act discretely and simplicially
on a coarse $PD(3)$ space.

\begin{remark}
The complex $\Si$ -- and hence $BS(p,q)$ -- can be uniformly properly embedded in a 
coarse $PD(3)$ space homeomorphic to $\R^3$.   To see this
we proceed as follows.  First take a proper PL embedding
$T\ra \R^2$ of the Bass-Serre tree into $\R^2$.   For each 
co-oriented edge $\oa{e}$ of $T\subset\R^2$  we take product
cell structure on the half-slab 
$P(\oa{e})\defeq \pi^{-1}(e)\times\R_+$ where $\R_+$ is given the usual cell structure.  
We now perform two types of gluings.   First, for each co-oriented edge
$\oa{e}$ we glue the half-slab $P(\oa{e})$ to $\Si$ by identifying
$\pi^{-1}(e)\times 0$ with $\pi^{-1}(e)\subset \Si$.   Now, for each
pair $\oa{e_1},\,\oa{e_2}$ of adjacent co-oriented edges, we glue 
$P(\oa{e_1})$ to $P(\oa{e_2})$ along $\pi^{-1}(v)\times\R_+$ where
$v=e_1\cap e_2$.   It is easy to see that after suitable subdivision
the resulting complex $X$ becomes a bounded geometry, uniformly
acyclic  $3$-dimensional PL manifold homeomorphic to $\R^3$.
\end{remark}

\medskip
\no {\bf Higher genus Baumslag-Solitar groups.} Note that $BS(p,q)$ is the
fundamental group of the following complex $K=K_1(p,q)$. Take the annulis $A$ 
with the boundary circles $C_1, C_2$. Let $B$ be another annulus with the 
boundary circles $C_1', C_2'$. 
Map $C_1', C_2'$ to $C_1, C_2$ 
by mappings $f_1, f_2$ of  degrees $p$ and $q$ respectively. Then $K$ 
is obtained by gluing $A$ and $B$ by $f_1\sqcup f_2$. Below we 
describe a ``higher genus'' generalization of this construction. Instead of 
the annulus $A$ take a surface $S$ of genus $g\ge 1$ with two boundary 
circles $C_1, C_2$. Then repeat the above construction of $K$ by gluing 
the annulus $B$ to $S$  via the mappings $C_1'\to C_1,
C_2'\to C_2$ of the degrees $p,q$ respectively. 
The fundamental group $G=G_g(p,q)$ of the resulting complex $K_g(p,q)$ 
has the presentation
$$
\< a_1,b_1,..., a_g, b_g, c_1, c_2, t: [a_1,b_1]...[a_g,b_g]c_1c_2=1, 
t c_2^q t^{-1}= c_1^p\>. 
$$
One can show that the group $G_g(p,q)$ is torsion-free and Gromov-hyperbolic 
\cite{KK}.  Note that the universal cover $\tilde{K}$ of the complex $K_g(p,q)$ 
does not fiber over the Bass-Serre tree $T$ of the HNN-decomposition of $G$. Nevertheless  
there is a properly embedded $c_1$-invariant subcomplex in  $\tilde{K}$ which ($c_1$-invariantly) 
fibers over $T$ with the fiber homeomorphic to $\R$. This allows one to repeat the 
arguments given above for the group $BS(p,q)$ 
and show that the group $G_g(p,q)$ cannot act simplicially freely on a 
coarse $PD(3)$ space (unless $p=\pm q$). 
However in \cite{KK} we show that $G_g(p,q)$ contains a finite index subgroup 
isomorphic to the fundamental group of a compact 3-manifold with boundary.

\medskip
\no
{\bf Groups with too many coarsely non-separating Poincare duality subgroups.}
By Corollary \ref{pdn-1peripheral}, if $G$ is of type $FP$,
and $G\acts X$ is a free simplicial action on a coarse $PD(n)$
space, then there are only finitely many conjugacy classes of
coarsely non-separating maximal $PD(n-1)$ subgroups in $G$.

We now construct an example of a  $2$-dimensional group of
type $FP$ which has infinitely many conjugacy classes of
coarsely non-separating maximal surface subgroups; this
example does not fit into any of the classes described above.
Let $S$ be a $2$-torus with one hole, and let $\{a,\,b\}\subset H_1(S)$ 
be a set of generators.  Consider a sequence of embedded loops $\ga_k\subset S$
which represent $a+kb\in H_1(S)$, for $k=0,1,\ldots$.  Let 
$\Si$ be a $2$-torus with two holes.   Glue the boundary torus of $S\times S^1$
homeomorphically to one of the boundary tori of $\Si\times S^1$
so that the resulting manifold $M$ is not Seifert fibered.
Consider the sequence $T_k\subset M$ of embedded 
incompressible tori corresponding
to $\ga_k\times S^1\subset S\times S^1\subset M$.   Let 
$L\subset\pi_1(M)$ be the infinite cyclic subgroup 
generated by the homotopy class of $\ga_0$.  Finally,
we let $G$ be the double of $\pi_1(M)$ over the cyclic
subgroup $L$, i.e. $G\defeq \pi_1(M)*_L \pi_1(M)$.
Then the reader may verify the following:

1.  Let $H_i\subset \pi_1(M)\subset G$ be the image of the
fundamental group of the torus $T_i$ for $i>0$ (which is well-defined
up to conjugacy).  Then each $H_i$ is maximal in $G$,
and the $H_i$'s are pairwise non-conjugate in $G$.

2. Each  $H_i\subset\pi_1(M)$ coarsely separates $\pi_1(M)$
into precisely two  deep  components. 

3. For each $i>0$, the subgroup $H_i\subset \pi_1(M)$ 
coarsely separates some conjugate of $L$ in $\pi_1(M)$.

4.  It follows from 3 that $H_i$ is coarsely non-separating
in $G$ for $i>0$. 

5. $G$ is of type $FP$ and has dimension $2$. 

\no
Therefore $G$ cannot act freely simplicially on a coarse $PD(3)$ space. 

\subsection{Appendix: Coarse Alexander duality in brief}

We will use terminology and notation from section
\ref{geomprelim}.

\begin{theorem}
\label{shortad}
Let $X$ and $Y$ be bounded geometry uniformly acyclic 
metric simplicial complexes,
where $X$ is an $n$-dimensional PL manifold.    Let $f:C_*(Y)\ra C_*(X)$
be a uniformly proper chain map, and let $K\subset X$ be the support
of $f(C_*(Y))\subset C_*(X)$.   For every $R$ we may compose the Alexander
duality isomorphism $A.D.$ with the induced map on compactly supported
cohomology:
\begin{equation}
\label{compa}
\t H_{n-k-1}(X-N_R(K))\stackrel{A.D.}{\lra}H^k_c(N_R(K))
\stackrel{H^k_c(f)}{\RA}H^k_c(Y);
\end{equation}
we call this composition $A_R$.  Then

1. For every $R$ there is an $R'$ so that
\begin{equation}
\label{mono}
Ker(A_{R'})\subset Ker(\t H_{n-k-1}(X-N_{R'}(K))\ra \t H_{n-k-1}(X-N_R(K))).
\end{equation}

2. $A_R$ is an epimorphism for all $R\geq 0$.

3. All deep components of $X-K$ are stable, and their number is
$1+rank(H^{n-1}_c(Y))$.

4.  If $Y$ is an $(n-1)$-dimensional manifold, then
for all $R$ there is a $D$ so that any point in $N_R(K)$ lies
within distance $D$ of both the deep components of 
$X-N_R(K)$.

The functions $R'=R'(R)$ and $D=D(R)$ depend only on the geometry
of $X$ and $Y$ (via their dimensions and acyclicity functions), and on
the coarse Lipschitz constant and distortion of $f$. 
\end{theorem}
\proof
{\em Step 1.}  We construct a coarse Lipschitz chain map $g:C_*(X)\ra C_*(Y)$ 
as follows. For each vertex $x\in X, y\in Y$ we let $[x], [y]$ denote the  
corresponding element of $C_0(X), C_0(Y)$. 
To define $g_0:C_0(X)\ra C_0(Y)$ we map $[x]$ for 
each vertex $x\in X\subset C_0(X)$
to $[y]$, where we choose a vertex $y\in Y\subset C_0(Y)$ 
for which the distance $d(x,Support(f(y)))$ is minimal, and 
extend this homomorphism $\Z$-linearly to a map $C_0(X)\to C_0(Y)$.  
Now assume inductively that $g_j:C_j(X)\ra C_j(Y)$ has been defined by 
$j<i$.   For each $i$-simplex $\si\in C_i(X)$, 
we define $g_i(\si)$ to be a chain 
bounded by $g_{i-1}(\D\si)$ (where $Support(g_i(\si))$ lies
inside the ball supplied by the acyclicity function of $Y$).
Using a similar inductive procedure to construct chain homotopies,
one verifies:

a)  For every $R$ there is an $R'$ so that the composition
\begin{equation}
\label{htoinclusion}
C_*(N_R(K))\stackrel{g_*}{\ra}C_*(Y)\ra C_*(K)\ra C_*(N_{R'}(K))
\end{equation}
is chain homotopic to the inclusion by an $R'$-Lipschitz chain
homotopy with displacement $<R'$.

b)  There is a $D$ so that 
$$
C_*(Y)\stackrel{f}{\ra}C_*(K)\stackrel{g}{\ra}C_*(Y)
$$
is a chain map with displacement at most $D$ and $g\circ f$
is chain homotopic to $id_{C_*(Y)}$ by a $D$-Lipschitz chain map
with displacement $<D$.

{\em Step 2.}  Pick $R$, and let $R'$ be as in a) above.  If
$$
\al\in Ker(H^k_c(N_{R'}(K))\stackrel{H^k_c(f)}{\RA}H^k_c(Y)), 
$$
then $\al$ is in the kernel of the composition
$$
H^k_c(N_{R'}(K))\stackrel{H^k_c(f)}{\RA}H^k_c(Y)\stackrel{H^k_c(g)}{\RA}
H^k_c(N_R(K))
$$
which coincides with the restriction $H^k_c(N_{R'}(K))\ra H^k_c(N_R(K))$
by a) above.   Similarly, the composition 
$$
H^k_c(Y)\stackrel{H^k_c(g)}{\RA}H^k_c(N_R(K))\stackrel{H^k_c(f)}{\RA}H^k_c(Y)
$$
is the identity, so $H^k_c(f)$ is an epimorphism.  Applying the Alexander
duality isomorphism to these two assertions we get parts 1 and 2.

{\em Step 3.}  Let $C$ be  a deep component of $X-K$.   Suppose
$C_1,\,C_2$ are deep components of  $X-N_R(K)$ with $C_i\subset C$.
Picking points $x_i\in C_i$, the difference $[x_1]-[x_2]$ determines
an element of $\t H_0(X-N_R(K))$ lying in 
$Ker(\t H_0(X-N_R(K))\ra \t H_0(X-K)$.  Hence 
$$
A_R([x_1]-[x_2])=A_0(p_R([x_1]-[x_2]))=A_0(0)=0
$$
where $p_R:\t H_0(X-N_R(K))\ra \t H_0(X-K)$ is the projection.
Since $C_1$ and $C_2$ are deep, for any $R'\geq R$ there is a
$c\in \t H_0(X-N_{R'}(K))$ which projects to $[x_1]-[x_2]\in \t H_0(X-N_R(K))$.
But then $A_{R'}(c)=0$ and part 1 forces $[x_1]-[x_2]=0$.
This proves that $C_1=C_2$, and hence that all deep components
of $X-K$ are stable.  The number of deep components of $X-K$
is 
$$
1+rank(\underset{\underset{R}{\longleftarrow}}{\lim}\t H_0(X-N_R(K)),
$$
and by part 1 this clearly coincides with $1+rank(H^{n-1}_c(Y))$.
Thus we have proved 2.

{\em Step 4.}
To prove part 4, we let $C_1,\,C_2$ be the two deep components
of $X-K$ guaranteed to exist by part 3.   Pick $x\in N_R(K)$, and
let $R'$ be as in part 1.    Since $f$ is coarse Lipschitz chain map,
there is a $y\in Y$ with $d(x,Support(f([y])))<D_1$ where
$D_1$ is independent of $x$ (but does depend on $R$).  
Choose a cocycle $\al\in C^{n-1}_c(Y)$
representing the generator of $H^{n-1}_c(Y)$ which is supported
in an $(n-1)$-simplex containing $y$.   Then the image $\al'$  of $\al$
under $C^{n-1}_c(Y)\stackrel{C^{n-1}_c(g)}{\lra}C^{n-1}_c(N_{R'}(K))$
is a cocycle supported in $B(x,D_2)\cap N_{R'}(K)$ where $D_2$ depends 
on $R'$ but is independent of $x$.  Applying the Alexander duality 
isomorphism\footnote{That is ultimately induced by taking the 
cap product with the fundamental class of $H^{lf}_n(X)$, the locally 
finite homology group of $X$.} to $[\al']\in H^{n-1}_c(N_{R'}(K))$, 
we get an element $c\in \t C_0(X-N_{R'}(K))$
which is supported in $B(x,D_2+1)\cap (X- N_{R'}(K))$,
and which maps under $A_{R'}$ to $[\al]\in H^{n-1}_c(Y)$.  Picking
$x_i\in C_i$ far from $K$, we have 
$[x_1]-[x_2]\in \t H_0(X-N_{R'}(K))$ and $A_{R'}([x_1]-[x_2])=\pm[\al]$.
By part 1 it follows that the images of $c$ and $[x_1]-[x_2]$
under the map $\t H_0(X-N_{R'}(K))\ra\t H_0(X-N_R(K))$ 
coincide up to sign.  In other words, $support(c)\cap C_i\neq\emptyset$,
so we've shown that $d(x,C_i)<D_2$ for each $i=1,2$.
\qed

\bibliography{refs}
\bibliographystyle{siam}
\addcontentsline{toc}{subsection}{References}

\medskip
\no

\parbox{3.5in}{Michael Kapovich:\\
~\\
Department of Mathematics\\
University of Utah\\
Salt Lake City, UT 84112-0090\\
kapovich@math.utah.edu}  \parbox{3.5in}{Bruce Kleiner:\\
~\\
Department of Mathematics\\
University of Michigan\\
 Ann Arbor, MI 48109\\
bkleiner@math.lsa.umich.edu}

\end{document}